\documentclass{amsart}

\newtheorem{Theo}{Theorem}[section]

\newtheorem{Lem}[Theo]{Lemma}

\newtheorem{Prop}[Theo]{Proposition}

\newtheorem{Cor}[Theo]{Corollary}

\newcommand{\sousfleche}[2]{
\renewcommand{\arraystretch}{0.5}\begin{array}[t]{c}#1\\{\scriptstyle#2}
\end{array}\renewcommand{\arraystretch}{1}}

\newcommand{\Hom}{\mathop{\mbox{\rm Hom}}\nolimits}

\newcommand{\Ker}{\mathop{\mbox{\rm Ker}}\nolimits}

\newcommand{\Gen}{\mathop{\mbox{\rm Gen}}\nolimits}

\newcommand{\Max}{\mathop{\mbox{\rm Max}}\nolimits}

\newcommand{\dist}{\mathop{\mbox{\rm dist}}\nolimits}

\newcommand{\myqed}{\ \hfill \rule{2mm}{2mm}}

\begin{document}  
\title[The Stable Rank and Swan's Problem]{The Stable Rank of 
Topological Algebras and a Problem of R.G. Swan}
\author{C. BADEA}
\address{Universit\'e de Lille I, F--59655 Villeneuve d'Ascq, France}
\email{badea@gat.univ-lille1.fr}


\begin{abstract}
Various notions of stable ranks are studied for topological algebras. Some
partial answers to R.G. Swan's problem (Have two Banach or good Fr\'echet algebras
as in the density theorem in $K$-theory the same stable rank~?) are obtained. For
example, a Fr\'echet dense $\ast$-subalgebra $A$ of a $C^\ast$-algebra $B$, closed
under $C^\infty$-functional calculus of self-adjoint elements, has the same
Bass stable rank as $B$. 
\end{abstract}



\maketitle
\tableofcontents
\section{Introduction}

\subsection{Preamble} Over the last ten years, there has been some interest in
studying the stable rank of Banach algebras. For the algebra of all continuous
complex or real-valued functions on a compact space many stable
ranks are related to the dimension of the space. Thus the theory of
stable ranks of $C^\ast$-algebras can be viewed as a kind of
noncommutative theory of dimension. 

We study in this paper various notions of stable ranks. The starting point was an
open problem of Richard G. Swan \cite[p.206]{Swa77} about the coincidence of the
Bass stable rank $Bsr$ (and other notions of stable ranks) of some topological
algebras (or rings) known to have the same $K$-theory by the density theorem. 

In the present paper (cf. Section \ref{main result}) we prove, 
among other things,
the following result~:

\begin{Theo}
Let $A$ be a dense and full (i.e. spectral invariant) $\ast$-algebra of a
$C^{\ast}$-algebra $B$, which is a (unital) Fr\'echet $Q$-algebra in its own
topology and closed under $C^{\infty}$-functional calculus of self-adjoint
elements. Then $Bsr(A) = Bsr(B)$. 
\end{Theo}

\smallskip

This kind
of ''smooth subalgebras`` $A$ of $C^{\ast}$-algebras was considered recently
from many points of view by several authors ; cf. for instance \cite{Co96},
\cite{BlCu91} and  their references.
The closure under holomorphic functional calculus of a dense subalgebra in 
the $C^{\ast}$-algebra completion (spectral invariance) as well as closure under
$C^{\infty}$-functional calculus of self-adjoint elements are two requirements for
any reasonable theory of ''smooth subalgebras`` of $C^{\ast}$-algebras
\cite{BlCu91}. The above result shows that smooth subalgebras must have the same
Bass stable rank as their $C^{\ast}$-completion. The commutative example we will
always keep in mind is that of $C^\infty (X)$ and $C(X)$.

We refer to the table of contents for the organization of the paper. Some of the
results of the present paper were announced in \cite{Bad94}. We refer however to
the present paper for the definitive statements. 
%

\subsection{Stable ranks} The first notion of stable rank for rings was introduced
by H. Bass \cite{Bas64} in his work in
algebraic 
$K$-theory. In this paper we will work in the framework of complex topological
algebras. A topological algebra is a complex topological vector space 
endowed with an 
algebra structure for which the multiplication is jointly continuous.
A unital topological algebra $A$ is called a $Q$-{\it{algebra}\/} if 
its set
of invertibles $A^{-1}$ is open. $A$ is called {\it{good}\/} if it is a 
$Q$-algebra and $x\mapsto x^{-1}$ is continuous on $A^{-1}$. A {\it{Fr\'echet
algebra}\/} is a (locally convex) Fr\'echet space with an algebra structure for
which the multiplication is jointly continuous. Note that we do not
require Fr\'echet algebras to be $m$-convex, i.e. admitting a family of
submultiplicative seminorms giving the topology. Fr\'echet algebras
are not necessarily $Q$-algebras, while Banach
algebras are. By \cite[p. 258]{Hel93}, unital Fr\'echet $m$-convex
$Q$-algebras are good Fr\'echet algebras. However, a {\it{normed
algebra}\/} is a normed space with an algebra structure and the norm is supposed
to be submultiplicative. 

By definition, the {\it{Bass stable rank}\/} $Bsr(A)$ of 
the algebra $A$
is the smallest positive integer $n$ such that the following condition

\begin{quote}
$(Bsr)_n$ \quad for every $(a_1, \ldots ,a_{n+1})\in Lg_{n+1}(A)$,
there exists
$(c_1,\ldots ,c_n) \in A^n$ such that $(a_1+c_1a_{n+1},\ldots
,a_n+c_na_{n+1})\in Lg_n(A)$, 
\end{quote}

\noindent holds, or infinity if no such number $n$ exists. 
Here $Lg_n(A)$ denotes the
set of $n$-tuples $(a_1,\ldots,a_n)$ of $A^n$ generating $A$ as a
left ideal :
$Aa_1+\cdots + Aa_n = A$. Note that in some of the earlier papers
the notion of the Bass stable rank was defined using different
indexing conventions. The condition
$(Bsr)_n$ was devised by Bass in order to determine values of $n$ for which every
matrix in  $GL_n(A)$ can be   row reduced by addition operations with coefficients
from $A$ to a matrix with the same last row and column as the identity matrix and to
obtain stability results in $K$-theory. Here $GL_n(A)$ is the group of invertible
elements of $M_n(A)$, the set of all $n\times n$ matrices with entries in $A$. It can
be proved that $(Bsr)_n$ implies $(Bsr)_{n+1}$ (\cite{Vas71},\cite{Kru73}) and that
the Bass stable rank of $A$ equals the stable rank of the opposed algebra $A^\circ$
(\cite{Vas71},\cite{War80}). That is, the left Bass stable rank defined
above is equal to the right Bass stable rank which can be defined in a similar way
using the set $Rg_n(A)$ of all $n$-tuples $(a_1,\ldots ,a_n)$
generating $A$ as a right ideal.

Inspired by one of the definition of the covering dimension of a topological space
$X$ in terms of the $C^\ast$-algebra $C(X)$ of all continuous functions on it,
Rieffel \cite{Rie83} introduced the notion of topological stable rank as follows~:
the {\it{topological stable rank}\/} $tsr(A)$ of the topological algebra $A$ is the
smallest positive integer $n$ such that $Lg_n(A)$ is dense in $A^n$,
or infinity if no such number exists. Here the topology of $A^n$ is
the product topology. If
$Lg_n(A)$ is dense in $A^n$, then $Lg_m(A)$ is dense in $A^m$ for every $m\geq n$.
The right topological stable rank $rtsr$ can be defined by considering the set
$Rg_n(A)$ instead of $Lg_n(A)$. We use here $tsr(A)$ instead of $\ell tsr(A)$ for
the left topological stable rank. The left and the right topological stable ranks
coincide for Banach algebras with a continuous involution. It is an open question
\cite[Question 1.5]{Rie83} if $tsr(A)$ equals $rtsr(A)$ for all Banach algebras $A$.
We refer to \cite{Rie83} for several properties of $tsr$. In particular, we have 
$Bsr(A) \leq tsr(A)$ for all Banach algebras $A$, with
equality for $C^\ast$-algebras as was shown by Herman and Vaserstein \cite{HeVa84}.
For the unital commutative $C^\ast$-algebra $C(X)$ we have
$$Bsr(C(X)) = tsr(C(X)) = [(\dim X)/2] + 1\;,$$
where $\dim X$ is the \v{C}ech-Lebesgue covering dimension of the space $X$
\cite{Pea75}.

Now we mention briefly other notions of stable ranks. The {\it{connected stable
rank}\/} $csr(A)$ of the Banach algebra $A$ is \cite{Rie83} the least integer $n$ such
that $GL_k(A)_0$ acts transitively (by left multiplication) on $Lg_k(A)$ for 
every  $k\geq n$, or, equivalently \cite{Rie83}, the least integer $n$ such that
$Lg_k(A)$ is connected for every $k\geq n$. We put $csr(A) =
\infty$ if no such $n$ exists. This notion is left-right symmetric \cite{CoLa86b}.
Here $GL_k(A)_0$ is the connected component of $GL_k(A)$ containing the identity.
The number $n$ belongs to the {\it{absolute connected stable range}\/} of $A$ if
for every nonempty, connected open subset $V$ of $A^n$, $V\cap Lg_n(A)$ 
is nonempty
and connected. The {\it{absolute connected stable rank}\/} $acsr(A)$ is
\cite{Nis86} the least integer $n$ such that,
for every $m\geq n$, $m$ belongs to
the absolute connected stable range of $A$ ; if no such integer exists
we set $acsr(A) = \infty$. Nistor \cite{Nis86} proved that
$acsr(A) = Bsr (C(I)\otimes A)$ ($= tsr(C(I)\otimes A)$), for all
$C^\ast$-algebras $A$, where $I = [0,1]$. A generalization of this
formula will be given in section 3.

The {\it{left (right) general stable rank}\/} of $A$ is defined 
\cite{Rie83} as
the smallest integer $n$ such that $GL_k(A)$ acts on the left 
(right) transitively
on $Lg_k(A)$ for all $k\geq n$, or, equivalently (for $C^\ast$-algebras), that the
cancellation property holds for finitely generated projective $A$-modules. For
instance, the right general stable rank is the smallest integer $n$ such that
$W\oplus A \cong A^k$ for some $k\geq n$ implies $W\cong A^{k-1}$, whenever $W$ is a
finitely generated projective left $A$-module. This notion is left-right symmetric
\cite{CoLa86b} and we will denote by $gsr(A)$ the common value. Another
characterization (for $C^\ast$-algebras) in terms of the equality of $Lg_k(A)$ with
the space $Lc_k(A)$ of last columns of the various matrices in $GL_k(A)$ was given
in \cite{Rie87}. If no such integer exists we set $gsr(A) = \infty$.
We have
$$gsr(A) \leq csr(A) \leq 1 + Bsr(A) \leq 1 + tsr(A)\;.$$
%

\subsection{Swan's problem} The Bass and topological stable ranks are useful for
stability results in the $K$-theory of topological algebras. We state here some of
them in terms of $Bsr$ ; since $Bsr(A) \leq tsr(A)$, we can replace in these
statements the Bass stable rank with the topological one. It was proved by G. Corach
and A.R. Larotonda \cite{CoLa86a} that, for a Banach algebra $A$ with unit $e$,
the map
$$
GL_n(A)/GL_n(A)_0\ni a\mapsto 
\left(
\begin{array}{cc} 
a & 0\\ 
0 & e
\end{array}
\right)
\in GL_{n+1}(A)/GL_{n+1}(A)_0$$
is bijective for $n\geq 1 + Bsr(A)$. Thus the topological
$K_1$-group $K_1(A)$ of $A$, defined as the direct limit of
$GL_n(A)/GL_n(A)_0$ under these inclusions, stabilizes if $Bsr(A)$
is finite. A similar result can be stated for the (topological)
$K_0$ group of a Banach algebra. Indeed, $K_0(A)$ can be written
\cite{Kar79} as the direct limit of $\bar P_{2n}(A)$ under the inclusions
$$p\mapsto
\left(
\begin{array}{ccc}
p & 0 & 0\\
0 & e & 0\\
0 & 0 & 0
\end{array}
\right) \; ,
$$
where $\bar P_{2n}(A)$ is the quotient of $P_{2n}(A)$, the subset of indempotents
of $M_{2n}(A)$, given by the equivalence relation
$$p_1 \approx p_2 \Leftrightarrow \exists \alpha\in GL_{2n}(A),
\quad \alpha p_1\alpha^{-1} = p_2\;.$$
Then the map $\bar P_{2n}(A)\rightarrow\bar P_{2n+2}(A)$ is bijective for 
$n\geq Bsr(A)$ \cite{Cor86}.
 
Consider a homomorphism of unital Banach (or Fr\'echet)
algebras $f : A\rightarrow B$ having a dense
%
range. The isomorphism result
known in $K$-theory as the density theorem may be stated \cite{Bos90} for good
Fr\'echet algebras.
\begin{Theo}[Density Theorem.] For a continuous dense and full morphism $f :
A\rightarrow B$ of two good Fr\'echet algebras $A$ and $B$, the maps $f_\ast :
K_i(A)\rightarrow K_i(B)$, $i = 0,1$, are isomorphisms.
\end{Theo}
About the terminology : A morphism $f : A\rightarrow B$ is called {\it{dense}}
if it has a dense image and {\it{full}} (spectral invariant) if $a$ is invertible
in $A$ whenever  $f(a)$ is invertible in $B$. When $f$ is the inclusion map of a
subalgebra $A$ of the algebra $B$, the above conditions reduce to the fact that $A$
is a dense and full (=``pleine'') subalgebra of $B$ \cite{Bou67}, or that $A$ and
$B$ form a Wiener pair according to Naimark's \cite{Nai72} terminology. A
subalgebra $A$ is {\it{closed under the holomorphic functional calculus}\/} in $B$ 
\cite{Co81} if,
for every $a\in A$ and $g$ holomorphic in a neighborhood of the
spectrum of $a$ in
$B$, the element $g(a)$ of $B$ lies in $A$. Subalgebras closed 
under the holomorphic
functional calculus are full. Full Fr\'echet subalgebras of 
Fr\'echet $Q$-algebras
are closed under the holomorphic functional calculus \cite{Sch92}.
There are other names for full subalgebras ; we refer to
\cite{Bad94} and the references cited therein.
%

The question whether two good Fr\'echet algebras $A$ and $B$ as in the
density theorem have the same stable rank (at least for a suitable
definiton of the stable rank) arises. It was Richard G. Swan
\cite[p.206]{Swa77} who raised this question for the Bass stable
rank and for the projective stable rank (see \cite{Swa77} for a
precise definition) in a less restrictive setting of some
topological rings. In the present paper, Swan's problem means this
problem for some topological algebras. A more general question to
ask is when do two Banach (or good Fr\'echet) algebras have the same
$K$-theory and their $K$-groups stabilize at the same level (i.e.
they have the same stable rank). Note that one can have one without
another~: $A$ and $C(X(A))$ have the same $K$-theory for all
commutative Banach algebras $A$ \cite{Tay75}, but for $A$ the disc
algebra or the algebra $H^\infty({\bf D})$ of bounded analytic
functions on the unit disc, the algebras $A$ and $C(X(A))$ have
different stable ranks $Bsr$ (but same $tsr$ ) :
see \cite{JMW86}, \cite{Tre92}, \cite{Sua94}.
Here $X(A)$ denotes the maximal ideal space of the commutative Banach
algebra $A$.

We refer to the introduction of \cite{Bad94} for the motivation of Swan's
problem and for some known results. The main results of the present paper are 
proved in section \ref{main result}.
%
\section{Characterizing $Bsr$ and $tsr$}
A subalgebra of a unital subalgebra is assumed to be a
unital subalgebra, i.e. it contains the identity
element $e$ of  the larger algebra. A morphism between
two unital algebras is supposed to satisfy $f(e) = e$.

\subsection{The Bass and the bilateral Bass stable
ranks} We introduce the following modification of the
Bass stable rank. We call {\it{bilateral Bass stable rank}\/} $bBsr(A)$ 
of the unital
algebra $A$ the smallest positive integer $n$ such that the following
condition holds for $k\geq n$
\begin{quote}
$(bBsr)_k$ \quad for every $(a_1,\ldots ,a_{k+1})\in
Lg_{k+1}(A)$, there
exist $(c_1,\ldots ,c_k)\in A^k$, $(d_1,\ldots
,d_k)\in A^k$, such that
$(a_1+c_1a_{k+1}d_1, \ldots , a_k+c_ka_{k+1}d_k) \in Lg_k(A)$,
\end{quote}

\noindent or infinity, if no such number exists. We have 
$bBsr(A)\leq Bsr(A)$, with equality for commutative
algebras. It will be shown later on that there are
$C^\ast$-algebras $A$ with $bBsr(A)\neq
Bsr(A)$. It would be interesting to study which
properties of $Bsr$ are preserved by the bilateral
stable rank~; we don't know for instance if condition
$(bBsr)_k$ implies $(bBsr)_{k+1}$.

We start off with some characterizations for $Bsr$ and
$bBsr$. While for $bBsr$ we can characterize it
in terms of onto algebra morphisms, for $Bsr$ we need
to work with left $A$-module morphisms. For a left
$A$-module $M$ and integer $n$ we let $\Gen_n(M)$
denote the $n$-{\it{generator space}\/} of $M$, that
is the set of elements $(v_1,\ldots,v_n)$ of $M^n$
such that $(v_1,\ldots,v_n)$ generates (algebraically)
$M$ : $Av_1+\cdots + Av_n = M$. If ${}_AA$ denotes $A$
itself viewed as a left $A$-module, then
$\Gen_n({}_AA) = Lg_n(A)$. If we identify $M^n$ with
$\Hom_A(A^n,M)$ in the canonical way, then
$\Gen_n({}_AM)$ corresponds to the surjections.                              
%
We will need the following auxiliary result.

\begin{Lem}\label{Lemma 2.1} 
Let $A$ be a unital topological $Q$-algebra. Then $Lg_n(A)$ is
open in
$A^n$.
\end{Lem}

{\bf Proof :} Consider a neighborhood $V$ of $0$ in $A$ such that 
$e+V\subseteq
A^{-1}$ ($A^{-1}$ is open). Let ${{\mathbf{u}}} = (u_1,\ldots ,u_n) \in
Lg_n(A)$, with $w_1u_1 + \cdots + w_nu_n = e$ for certain elements
$w_1,\ldots, w_n$ in
$A$. We can choose a
neighborhood $W$ of $0$ in $A$ such that $w_1W + \cdots + w_nW\subseteq V$. We show
that every ${{\mathbf{v}}} = (v_1,\ldots ,v_n)$ with
${{\mathbf{v}}}-{{\mathbf{u}}}\in W^n$ belongs to $Lg_n(A)$. Indeed, for all
$i\in\{1,\ldots ,n\}$, we have
$$w_iv_i = w_iu_i + w_i(v_i-u_i)\in w_iu_i + w_iW\;.$$
Therefore 
$$w_1v_1 + \cdots + w_nv_n\in (e+w_1W + \cdots + w_nW)\subseteq e+V\;.$$
This implies that $\sum^n_{i=1}w_iv_i$ is invertible in $A$ and therefore
$(v_1,\ldots ,v_n)\in Lg_n(A)$. \myqed

\smallskip

The following characterization shows that the condition 
``$Bsr(A)\leq n$'' is 
equivalent to a certain lifting condition.

\begin{Theo}\label{Theorem 2.2} 
Let $A$ be a unital topological $Q$-algebra. The following two
assertions are equivalent : 
\begin{itemize}
\item[(i)] The condition $(Bsr)_n$ holds, that is $Bsr(A) \leq n$ ;
\item[(ii)] Every onto module morphism of left $A$-modules 
$f : {}_AA\rightarrow
{}_AM$ induces an onto mapping $f_n : Lg_n(A)\rightarrow \Gen_n(M)$.
\end{itemize}
\end{Theo}

{\bf Proof :} ``$\Rightarrow$''. Suppose that (i) holds and let 
$f : A\rightarrow M$
be a mapping as in (ii). Let $f_n({\mathbf{a}}) = (f(a_1),\ldots
,f(a_n))\in
\Gen_n(M)$, ${\mathbf{a}}\in A^n$. We get $\sum^n_{i=1}y_if(a_i) = f(e)$ for
certain $y_i\in A$, and thus $f(\sum^n_{i=1}y_ia_i-e) = 0$ since $f$ is a
module morphism. Denote $u = \sum^n_{i=1}y_ia_i-e$. Then $u$ is an element
of the left ideal $\Ker (f)$ of $A$ and also $(a_1,\ldots ,a_n,u)\in
Lg_{n+1}(A)$. By (i), there are elements $c_i$ in $A$, $1\leq i\leq n$, such
that ${\mathbf{x}} = (a_1 + c_1u,\ldots ,a_n + c_nu)\in Lg_n(A)$. Since
$u\in\Ker (f)$, we get $f_n({\mathbf{x}}) = f_n({\mathbf{a}})$ and thus $f_n :
Lg_n(A)\rightarrow \Gen_n(M)$ is onto.

``$\Leftarrow$''. Suppose that (ii) holds and consider $(a_1,\ldots
,a_{n+1})\in Lg_{n+1}(A)$. Let $J$ be the closed left ideal generated by
$a_{n+1}$ and let $M$ be the left $A$-module $A/J$. Then the canonical mapping
$\pi : A\rightarrow M$ is onto. By (ii) we obtain another onto mapping $\pi_n :
Lg_n(A)\rightarrow \Gen_n(M)$. Since $(a_1,\ldots ,a_{n+1})\in
Lg_{n+1}(A)$, we have
$x_1a_1+\cdots +x_na_n + x_{n+1}a_{n+1} = e$ for certain elements
$x_1,\ldots ,x_{n+1}$ in $A$. Therefore $x_1\pi(a_1) + \cdots +
x_n\pi(a_n) = \pi(e)$. Now, if $v\in M$, $v = \pi(a)$, $a\in A$, then
$$v = \pi (a) = a\pi (e) = ax_1\pi (a_1) + \cdots + ax_n\pi (a_n)\;.$$
Therefore $(\pi (a_1),\ldots ,\pi (a_n))\in \Gen_n(M)$ and, by the
surjectivity of 
$\pi_n$, we find ${\mathbf{x}} = (x_1,\ldots ,x_n)\in Lg_n(A)$ with
${\mathbf{x}}-{\mathbf{a}}\in J$, where ${\mathbf{a}} = (a_1,\ldots ,a_n)$.
Using Lemma \ref{Lemma 2.1}, there exists a neighborhood
$U$ of $0$ in $A^n$ such that
${\mathbf{a}}^{\ast}\in {\mathbf{x}}+U$ implies ${\mathbf{a}}^{\ast}\in
Lg_n(A)$. Since ${\mathbf{x}}-{\mathbf{a}}\in J$, we find the existence
of the
$n$-tuple ${\mathbf{c}}$ with elements of $A$ such that 
${\mathbf{x}}-{\mathbf{a}}
- {\mathbf{c}}a_{n+1}\in U$. It follows that
${\mathbf{a}}+{\mathbf{c}} a_{n+1}\in Lg_n(A)$, yielding $Bsr(A)\leq n$.
\myqed %

\smallskip

For the bilateral Bass stable rank, we need only algebra morphisms.

\begin{Theo}\label{Theorem 2.3} Let $A$ be a unital topological $Q$-algebra. 
The
following two assertions are equivalent :
\begin{itemize}
\item[(i)] The condition $(bBsr)_n$ holds ; 
\item[(ii)] Every onto unital algebra morphism $f : A\rightarrow B$, $B$ a
topological algebra, induces an onto mapping $f_n : Lg_n(A)\rightarrow Lg_n(B)$.
\end{itemize}
\end{Theo}

{\bf Proof :} The proof is similar to the proof of the previous theorem. For the
implication (i) $\Rightarrow$ (ii) we use the fact that $f$ is an algebra morphism,
instead of a module morphism. For the reverse implication we have to consider the
closed two-sided ideal generated by $a_{n+1}$, instead of the left ideal $J$.
 \myqed

\smallskip

In the case when the topological $Q$-algebra $A$ is commutative, we have $bBsr(A) =
Bsr(A)$.

\begin{Cor}\label{Corollary 2.4} 
Let $A$ be a unital commutative topological $Q$-algebra.
Then $Bsr(A)\leq n$ if and only if for every commutative unital topological 
algebra
$B$ and a unital onto algebra morphism $f : A\rightarrow B$, the induced mapping
$f_n : Lg_n(A)\rightarrow Lg_n(B)$ is also onto.
\end{Cor}
This equivalence has already been proved in \cite{EsOh67} for
commutative rings and in \cite{CoLa84} for commutative Banach
algebras.
%
\subsection{The topological stable rank}
For the topological stable rank we have the following characterization.

\begin{Theo}\label{Theorem 2.5}
Let $A$ be a unital good topological algebra.
The following three assertions are equivalent : 
\begin{itemize}
\item[(i)] $tsr(A) \leq n$ ;
\item[(ii)] for every onto $A$-module morphism $f : {}_AA \rightarrow {}_AM$
of left $A$-modules, the induced map $f_n : Lg_n(A)\rightarrow \Gen_n(M)$ has the
following property~: If ${\mathbf{a}} = (a_1,\ldots ,a_n)\in A^n$
satisfies
$f_n({\mathbf{a}})\in \Gen_n(M)$, then, for every open neighborhood $V$ of
$0$ in $A^n$, there exists ${\mathbf{a}}'\in Lg_n(A)$ such that
${\mathbf{a}}'-{\mathbf{ a}}\in V$ and $f_n({\mathbf{a}}') =
f_n({\mathbf{a}})$.
\item[(iii)] for every onto morphism $f : A\rightarrow B$ of topological algebras,
the induced map $f_n : Lg_n(A)\rightarrow Lg_n(B)$ has the following property : If
${\mathbf{a}} = (a_1,\ldots ,a_n)\in A^n$ satisfies $f_n({\mathbf{a}})\in
Lg_n(B)$, then, for every open neighbourhood $V$ of $0$ in $A^n$,
there exists
${\mathbf{a}}'\in Lg_n(A)$ such that ${\mathbf{a}}'-{\mathbf{a}}\in V$ and
$f_n({\mathbf{a}}') = f_n({\mathbf{a}})$.
\end{itemize}
\end{Theo}

{\bf Proof :} In what follows we will consider an element ${\mathbf{a}}$ of
$A^n$ both as a line and a row vector. We consider $M_n(A)$ as a topological
algebra with topology induced by that of $A$ \cite{Swa77}.

\smallskip

(i) $\Rightarrow$ (ii). Suppose $tsr(A) \leq n$. Denote by 
$f^n : M_n(A)\rightarrow
M_n(M)$ the canonical extension of the mapping $f$ over the matrix algebra with 
entries in $A$ onto the $M_n(A)$-module $M_m(M)$. Since $f$ is onto, $f^n$ is also
onto. Let ${\mathbf{a}}$ be as in (ii) and let $V$
%
be an arbitrary neighborhood of $0$ in $A$. We will denote again by $V$ the
corresponding neighborhood of $0$ in $M_n(A)$ or the neighborhood of 
$0$
in $A^n$. By \cite[Lemma 2.1]{Swa77}, if $A$ is good so is $M_n(A)$.

Let $V'$ be another neighborhood of $0$ (in $A^n$) such that 
$V'({\mathbf{a}} + V')
+ V'\subset V$. It exists an open neighborhood $U$ of $0$ in $A$ such that if 
$x\in M_n(A)$ has all entries in $U$, then $e+x$ is invertible and $(e+x)^{-1}\in
e+V'$. Set
$$W = \{ (e+x)f_n(a) : x\in M_n(A)\quad\hbox{\rm has all entries in}\; U\} .$$
Then $W\subseteq \Gen_n(M)$ since $e+x$ is invertible in $M_n(A)$. Thus 
${\mathbf{
a}}\in f^{-1}_n(W)$. Also $f^{-1}_n(W)$ contains an open neighborhood of 
${\mathbf{
a}}$ (this is clear if $f$ is continuous). Indeed, because
$f_n({\mathbf{a}})\in \Gen_n(M)$, there exists $(s_1,\ldots ,s_n)\in A^n$
such that  %
$$s_1f(a_1) + \cdots + s_n f(a_n) = f(e)\;.$$
Let $T$ be a neighborhood of $0$ in $A$ such that $Ts_i\subset U$ for all  $i$,
$1\leq i\leq n$. Let ${\mathbf{z}} = (z_1,\ldots , z_n)\in A^n$ be such
that 
$z_i-a_i\in T$ ($i=1,\ldots ,n$) and set ${\mathbf{y}} = (y_1,\ldots
,y_n) = f({\mathbf{z}})\in M^n$. We show that there exists $x\in
M_n(A)$, having all entries in $U$, such that ${\mathbf{y}} = (e+x)
f_n({\mathbf{a}})$. Indeed, for every $k\in\{ 1,\ldots ,n\}$, we have
$$f(z_k)-f(a_k) = (z_k-a_k) f(e) = \sum^n_{i=1}(z_k-a_k)s_if(a_i)$$
and $(z_k-a_k)s_i\in Ts_i\subset U$ for all $i$. Then
%
$$y_k = f(z_k) = \sum^n_{i=1}[\delta_{ik} + (z_k-a_k)s_i]f(a_i)$$          
where $\delta_{ik}$ is $1$ if $i = k$ and $0$ if not. 
Thus ${\mathbf{y}} = (e+x)f_n({\mathbf{
a}})$ and every entry of the matrix $x = [(z_k-a_k)s_i]_{ik}$ belongs to $U$.

The set $Lg_n(A)$ is dense in $A^n$, so it meets $f^{-1}_n(W)$. 
Choose ${\mathbf{
a}}''\in Lg_n(A)\cap f^{-1}_n(W)$ such that ${\mathbf{a}}''-{\mathbf{a}}\in
V'$. We have $f_n({\mathbf{a}}'') = (e+x) f_n({\mathbf{a}})$ for a suitable
$x\in M_n(A)$ with all entries in $U$. Consider ${\mathbf{a}}' =
(e+x)^{-1}{\mathbf{a}}''\in Lg_n(A)$. Then 
\begin{eqnarray*}
{\mathbf{a}}' -{\mathbf{a}} &=& {\mathbf{a}}'-{\mathbf{a}}''+ {\mathbf{
a}}''-{\mathbf{a}}\\
&=& \{(e+x)^{-1}-e\} {\mathbf{a}}'' + ({\mathbf{a}}''-{\mathbf{a}})\\
&\in & V'{\mathbf{a}}''+V'\\
&\subset & V'({\mathbf{a}}+V')+V'\subset V
\end{eqnarray*}

\noindent and
$$f_n({\mathbf{a}}') = f_n\bigl ((e+x)^{-1}{\mathbf{a}}''\bigr) =
(e+x)^{-1}f_n({\mathbf{a}}'') = f_n({\mathbf{a}})\;.$$

(ii) $\Rightarrow$ (iii). Let $B$ be a unital topological algebra and consider an
onto, continuous and unital algebra morphism $f:A\rightarrow B$. Then $B$ 
becomes a
left $A$-module with respect to the multiplication $a\ast b = f(a)b$, where $b\in
B$, $a\in A$. Indeed,
\begin{eqnarray*}
a\ast (b_1+b_2) &=& f(a)(b_1+b_2) = f(a)b_1+f(a)b_2\\
&=& a\ast b_1 + a\ast b_2 \; ,
\end{eqnarray*}
%
$(a_1+a_2)\ast b = a_1\ast b + a_2\ast b$ in a similar fashion and also
$$(a_1a_2)\ast b = f(a_1a_2)b = f(a_1)f(a_2)b = a_1\ast(a_2\ast b)\;.$$
Here $a,a_1,a_2$ are in $A$ and $b,b_1,b_2$ are in $B$. The operation of
multiplication is well-defined and makes the map $f : {}_AA\rightarrow
{}_AB$ an (onto) left $A$-module morphism : $f(a_1a_2) = f(a_1)f(a_2) =
a_1\ast f(a_2)$. The element ${\mathbf{b}} = (b_1,\ldots ,b_n)$ is in
$Lg_n(B)$ if and only if $Bb_1 + \cdots + Bb_n = B$. Because $f$ is onto,
this is equivalent to $f(A)b_1 + \cdots + f(A)b_n = B$ and thus to
$f(Aa_1 + \cdots + Aa_n) = B$, where $b_i = f(a_i)$, $a_i\in A$, $1\leq
i\leq n$. On the other hand, ${\mathbf{b}} = f({\mathbf{a}})\in \Gen_n(B)$
if and only if $A\ast f(a_1)+\cdots + A\ast f(a_n) = B$, that is, if and
only if $f(Aa_1+\cdots +A a_n) = B$. Therefore $\Gen_n({}_AB) = Lg_n(B)$ and
(iii) follows from (ii).

(iii) $\Rightarrow$ (i). Let ${\mathbf{a}} = (a_1,\ldots , a_n)\in A^n$
and let $V$ be a neighborhood of $0$ (in $A^n$). Let $J$ be the closed
two-sided ideal in $A$ generated by $a_n-e$. If $\pi$ denotes the
canonical projection of $A$ onto the topological algebra $B = A/J$,
we have $\pi (a_n) = \pi(e)$. Therefore $(\pi (a_1),\ldots
,\pi(a_n))\in Lg_n(B)$. By (iii), there exists ${\mathbf{c}} =
(c_1,\ldots ,c_n)\in Lg_n(A)$ such that 
$\pi({\mathbf{c}}) = \pi({\mathbf{a}})$ and ${\mathbf{a}} -{\mathbf{c}}\in
V$. Therefore $Lg_n(A)$ is dense in $A^n$, implying $tsr(A)\leq n$.
\myqed     
%

\smallskip

These characterizations for the Bass and the topological stable 
ranks (Theorem \ref{Theorem
2.2} and Theorem \ref{Theorem 2.5}) imply the inequality $Bsr(A)\leq tsr(A)$
between the  (left) Bass stable rank and the left topological stable rank. Of
course, using  $Bsr(A) = Bsr(A^\circ)$, we obtain $Bsr(A)$ not greater than the
minimum of the left and the right topological stable ranks of $A$. This inequality
was obtained (for Banach algebras) in \cite[p.305]{Rie83} and 
\cite[Theorem 3]{CoLa84} by comparing first (\underbar{left}) $Bsr$
with
\underbar{right} $tsr$.

The following theorem provides a characterization for the topological stable rank
similar to Bass' definition of his stable rank as well as a bilateral version. It
implies again the inequality between the (left) Bass stable rank and the left
topological stable rank.

\begin{Theo}\label{Theorem 2.6} 
Let $A$ be a unital good topological algebra.
The following three assertions are equivalent
\begin{itemize}
\item[(i)] $tsr(A)\leq n$ ;
\item[(ii)] For every open neighborhood $V$ of $0$ in $A$ and every
$(a_1,\ldots ,a_{n+1})\in Lg_{n+1}(A)$, there exists $(c_1,\ldots
,c_n)\in A^n$ such that 
$(a_1+c_1a_{n+1},\ldots ,a_n+c_na_{n+1})\in Lg_n(A)$ and
$c_ia_{n+1}\in V$ for
$i=1,\ldots ,n$.
\item[(iii)] For every open neighborhood $V$ of $0$ in $A$ and every
$(a_1,\ldots ,a_{n+1})\in Lg_{n+1}(A)$, there exist $(c_1,\ldots
,c_n)$, $(d_1,\ldots ,d_n)\in A^n$ such that
$(a_1+c_1a_{n+1}d_1,\ldots ,a_n+c_na_{n+1}d_n)\in Lg_n(A)$ and 
$c_ja_{n+1}d_j\in V$ for $j = 1,\ldots ,n$.
\end{itemize}
\end{Theo}
%

{\bf Proof :} (i) $\Rightarrow$ (ii). Suppose $tsr(A)\leq n$ and let $V$ be an open
neighborhood of $0$. Consider $(a_1,\ldots ,a_{n+1})\in Lg_{n+1}(A)$
and let $J$ be the closed left ideal generated by $a_{n+1}$ in $A$.
Then $A/J$ is a left
$A$-module. Denoting $\pi : A\rightarrow A/J$ the canonical projection, we obtain,
as in the proof of Theorem \ref{Theorem 2.2}, that  
$(\pi(a_1),\ldots ,\pi(a_n))\in
\Gen_n(A/J)$. By  Theorem \ref{Theorem 2.5}, there exists ${\mathbf{t}} =
(t_1,\ldots ,t_n)\in Lg_n(A)$ such that $\pi(t_i) = \pi(a_i)$ and
$t_i-a_i\in\frac{1}{2}V$ for $i=1,\ldots ,n$. Since $Lg_n(A)$ is open
in $A^n$ (Lemma \ref{Lemma 2.1}), there exists an open neighborhood
$U$ of $0$ in $A^n$ such that $U\subset \frac{1}{2}V$ and
${\mathbf{t}}-{\mathbf{x}}\in U$ implies
${\mathbf{x}}\in Lg_n(A)$. We have $a_i-t_i\in J$ and thus we can find
$c_i\in A$, $i=1,\ldots,n$, such that
${\mathbf{a}}+{\mathbf{c}}a_{n+1}-{\mathbf{t}}\in U$. It follows that
${\mathbf{a}}+{\mathbf{c}}a_{n+1} \in Lg_n(A)$. Moreover, 
\begin{eqnarray*}
c_ia_{n+1} &=& a_i+c_ia_{n+1}-t_i+t_i-a_i\\
&\in & U + \frac{1}{2}V\subset V \; ,
\end{eqnarray*}
for all $i=1,\ldots ,n$.

\noindent (ii) $\Rightarrow$ (iii). Choose $d_i = e$, $i=1,\ldots
,n$.

\noindent (iii) $\Rightarrow$ (i). Let ${\mathbf{a}} = (a_1,\ldots
,a_n)\in A^n$. Then  $(a_1,\ldots,a_n,e)\in Lg_{n+1}(A)$. Using
(iii), for every open neighborhood $V$ of $0$ in $A$, there exists
%
${\mathbf{z}} = (z_1,\ldots ,z_n)\in A^n$ such that ${\mathbf{a}} +
{\mathbf{z}}
\in 
Lg_n(A)$ and $z_i\in V$ for any $i=1,\ldots ,n$. Hence
${\mathbf{a}}'\in Lg_n(A)$ and ${\mathbf{a}}-{\mathbf{a}}'\in V$, where
${\mathbf{a}}' = (a_1+z_1,\ldots ,a_n+z_n)$. Thus $tsr(A) \leq n$. 
\myqed 

\smallskip

Theorem \ref{Theorem 2.5} can be improved for commutative Banach
algebras.

\begin{Prop}\label{Proposition 2.7} 
Let $A$ be a commutative unital Banach algebra.
Then $tsr(A)\leq n$ if and only if for every multiplicative linear functional
$\gamma : A\rightarrow {\bf C}$, the induced map $\gamma_n$ has the following
property : If ${\mathbf{a}} = (a_1,\ldots ,a_n)\in A^n$,
$\gamma_n({\mathbf{a}})\in {\bf C}^n\setminus\{ 0\}$, then, for every
$\varepsilon > 0$, there exists
${\mathbf{a}}'\in Lg_n(A)$ such that $\|{\mathbf{a}}'-{\mathbf{a}}\| <
\varepsilon$ and $\gamma_n({\mathbf{a}}') = \gamma_n({\mathbf{a}})$.
\end{Prop}

{\bf Proof :} By Theorem \ref{Theorem 2.5} it is sufficient to prove the ``if''
part. Assume the condition above holds for every character $\gamma$.
We want to show that $Lg_n(A)$ is dense in $A^n$. Let ${\mathbf{a}} =
(a_1,\ldots ,a_n)\in A^n$ and let $\varepsilon > 0$. If ${\mathbf{a}}
= 0$, then the element ${\mathbf{c}} = (\varepsilon /2,0,\ldots ,0)$
belongs to
$Lg_n(A)$ and we have 
$\|{\mathbf{a}}-{\mathbf{c}}\| < \varepsilon$. Suppose now $a_k\neq 0$,
for a certain $k$ with $1\leq k\leq n$, and let $J$ be the ideal generated by
$a_k-2\| a_k\| e$. Since $\lambda = 2\| a_k\| >
\| a_k\|$, $\lambda$ is not in the spectrum of $a_k$. Therefore $J$
is a proper ideal. Let  $J_{\max}$ be the maximal ideal of $A$ containing $J$. Then
$A/J_{\max}$ is a Banach field and thus isomorphic with ${\bf C}$.
%
Let $\pi$ be the canonical projection of $A$ onto ${\bf C} = A/J_{\max}$. Then
$\pi$ is linear and multiplicative and 
$$\pi(a_k) = \pi(a_k-2\| a_k\| e+2\| a_k\| e) =
2\| a_k\| \pi(e)$$
We obtain $(\pi(a_1),\ldots ,\pi(a_n))\in{\bf C}^n\setminus \{0\}$.
Therefore there is ${\mathbf{c}} = (c_1,\ldots,c_n)\in Lg_n(A)$ such
that $\pi({\mathbf{c}}) =
\pi({\mathbf{a}})$ and $\| {\mathbf{a}}-{\mathbf{c}}\| < \varepsilon$. The
set $Lg_n(A)$ is thus dense in $A^n$, yielding $tsr(A) \leq n$. 
\myqed 
%

\subsection{Other stable ranks for normed algebras}
Based upon \cite{Rob80}, consider the {\it{unitary stable rank}\/} $usr(A)$ of
a  $C^\ast$-algebra $A$ which is, by definition, the least integer $n$ (or
infinity) with the property that for every $(a_1,\ldots, a_{n+1})\in
Lg_{n+1}(A)$ there exist $c_1,\ldots ,c_n$ in the unitary group of
$A$ such that
$(a_1+c_1a_{n+1},\ldots ,a_n+c_na_{n+1})\in Lg_n(A)$. We have
$Bsr(A)\leq usr(A)$. A.G. Robertson \cite{Rob80} has proved that,
for a $C^\ast$-algebra $A$ and for 
$n=1$, the condition $Bsr(A) = n$ is equivalent to $usr(A) = n$. It is an open
question \cite[Question 3.2]{Rie83} to extend Robertson's result in a useful way
for $n\geq 2$.

For a normed
algebra $A$ we
define the
{\it{norm
one stable rank}\/}
$nsr(A)$ as the
least integer $n$
such that, for every $(a_1,\ldots ,a_{n+1})\in Lg_{n+1}(A)$, there
are some elements
$c_1,\ldots ,c_n$ in $A$ of norm one such that
$(a_1+c_1a_{n+1},\ldots ,a_n+c_na_{n+1})
\in Lg_n(A)$, or infinity if no such number exists. We define the {\it{small
norm stable rank}\/} $ssr(A)$ as the least integer $n$ (or infinity) such that for every
$\varepsilon > 0$ and every $(a_1,\ldots ,a_{n+1})\in Lg_{n+1}(A)$,
there are some  elements $c_1,\ldots ,c_n$ in $A$ such that
$\| c_i\| < \varepsilon$ ($i=1,\ldots ,n$) and
$(a_1+c_1a_{n+1},\ldots ,a_n+c_na_{n+1})\in Lg_n(A)$.
%
\begin{Prop}\label{Proposition 2.8}
For every unital normed good algebra $A$ we have
$$bBsr(A)\leq Bsr(A)\leq tsr(A)\leq ssr(A)\leq nsr(A)\;.$$
For a $C^\ast$-algebra $A$ we have $nsr(A)\leq usr(A)$. Moreover, if one of 
the numbers $Bsr(A)$, $tsr(A)$, $ssr(A)$, $nsr(A)$, $usr(A)$ is equal to one, then
all these numbers are equal to one.
\end{Prop}

{\bf Proof :} Suppose $nsr(A)\leq n$. Let $(a_1,\ldots ,a_{n+1})\in
Lg_{n+1}(A)$ and  $\varepsilon > 0$. Then, for every positive
integer $k$, $(a_1,\ldots ,\frac{1}{k} a_{n+1})\in Lg_{n+1}(A)$. 
Therefore, there are $d_1,\ldots ,d_n$ in
$A$ of norm one such that $(a_1 + \frac{d_1}{k}a_{n+1},\ldots , a_n +
\frac{d_n}{k}a_{n+1})\in Lg_n(A)$. Choosing $k$ such that $k > 1/\varepsilon$, we
obtain $ssr(A) \leq n$. Choose $\varepsilon > 0$ and let
$(a_1,\ldots ,a_{n+1})\in Lg_{n+1}(A)$. Consider $\varepsilon' =
\varepsilon /\| a_{n+1}\|$. Since $ssr(A)\leq n$,
there exists $(c_1,\ldots ,c_n)\in A^n$ such that
$(a_1+c_1a_{n+1},\ldots ,a_n+c_na_{n+1})\in Lg_n(A)$ and $\|
c_i\| <
\varepsilon'$. Therefore $\| c_ia_{n+1}\| \leq \|
c_i\|\| a_{n+1}\| < \varepsilon$ for
$i=1,\ldots ,n$. By Theorem \ref{Theorem 2.6} we have $tsr(A)\leq n$.

The inequality $nsr(A) \leq usr(A)$ for $C^\ast$-algebras is clear. If one of the
values of the stable ranks mentioned in this proposition is equal to one, then
$Bsr(A) = 1$ and, by the above mentioned result due to Robertson, $usr(A) = 1$. 
Thus all stable ranks are one. 
\myqed 
%

\smallskip

We have $bBsr(A) = Bsr(A)$ for all commutative $A$. We also have equality for all
(normed) algebras with the following finiteness type property :
every two-sided ideal generated by $x$ (a principal ideal) is equal
to the left ideal generated by the same element. Theorem
\ref{Theorem 2.6} shows that the bilateral topological stable rank 
(defined by condition (iii) in Theorem \ref{Theorem 2.6}) coincides
with the usual $tsr$. For Bass stable rank this is not true.

\begin{Prop}\label{Proposition 2.9} For every unital simple $C^\ast$-algebra $A$ we have
$bBsr(A) = 1$. There are unital simple $C^\ast$-algebras satisfying $Bsr(A) = 1$ 
or $Bsr(A) = \infty$. There are non-unital simple $C^\ast$-algebras with 
$Bsr(A) = 2$.
\end{Prop}

{\bf Proof :} Since $A$ is simple, there are no proper closed two-sided ideals.
Therefore, the existence of onto morphisms $A\rightarrow A/J$, $J$ closed, implies
$J = \{0\}$ or $J = A$. In both cases the mapping $A^{-1}\rightarrow (A/J)^{-1}$ is
onto. By (the proof of) Theorem \ref{Theorem 2.3} we have $bBsr(A) = 1$.

There are several examples of simple $C^\ast$-algebras with $Bsr(A) = 1$~: see for
instance the references of \cite{BDR}. It follows from \cite[\S 2.2]{Cun77} and
\cite[\S 6.5]{Rie83} that a simple unital and purely infinite $C^\ast$-algebra $B$
has infinite Bass stable rank. Every such $B$
admits a full (nonunital but $\sigma$-unital) hereditary $C^\ast$-algebra of Bass
stable rank 2  (cf~. \cite{ElH}). A simple unital and purely
infinite $C^\ast$-algebra $B$ has infinite Bass stable rank ; however, it has real
rank $0$, (cf. \cite{Zha90}). We refer to the last section of the
present paper for the definition of the Bass stable rank 
for non-unital algebras. \myqed 

\smallskip

In particular, we can have $bBsr(A)\neq Bsr(A)$. The following result gives a
bilateral characterization of $Bsr$ for $C^\ast$-algebras.

\begin{Theo}\label{Theorem 2.10} 
Let $A$ be a unital $C^\ast$-algebra. The Bass
stable rank of $A$ is the least integer $n$ (or infinity) such that the following
condition
\begin{quote}
$(pbsr)_n$ \quad for every $(a_1,\ldots ,a_{n+1})\in Lg_{n+1}(A)$, there
exist $(c_1,\ldots ,c_n)\in A^n$ and $(d_1,\ldots ,d_n)\in A^n$, such that $d_i$
are positive, $1\leq i\leq n$, and $(a_1+c_1a_{n+1}d_1,\ldots
,a_n+c_na_{n+1}d_n)\in Lg_n(A)$  \end{quote}
holds.
\end{Theo}

{\bf Proof :} If $A$ satisfies $(Bsr)_n$, then $A$ satisfies $(pbsr)_n$ (just put
$d_i = e$, $1\leq i\leq n$).

We have $Bsr(A) = tsr(A)$ \cite{HeVa84}. We use the method from \cite{HeVa84} to
prove that $tsr(A)\leq n$ if condition $(pbsr)_n$ holds. We can suppose $n$ finite.
Let $\varepsilon$ be a given number, $0 <\varepsilon < 2$. Set $\varepsilon' =
\varepsilon^2/4 < 1$. Let $b_1,\ldots ,b_n$ be elements of $B$. Set
$b_0 =
\sum^n_{i=1}b^\ast_ib_i$ and $b_{n+1} = (e-b_0/\varepsilon')^{+}$, the positive
part of $e - b_0/\varepsilon'$. Then \cite[Lemma]{HeVa84} we have
$(b_1,\ldots ,b_{n+1})\in Lg_{n+1}(A)$. Therefore
$(b_1+c_1b_{n+1}d_1,\ldots ,b_n+c_nb_{n+1}d_n)\in Lg_n(A)$ for
suitable $c_i$ in $A$ and suitable positive elements $d_i$ in $A$,
$i=1,\ldots ,n$. Let $K \geq (1/\varepsilon')\max$ ($\|
c_i\| : 1 \leq i\leq n$) be fixed. Define $v = e+Kb_{n+1}d$,
where $d = d_1 + \ldots + d_n$ is positive. Then $v \geq e$, so  it
is invertible. Let $b'_i = (b_i+c_ib_{n+1}d_i)v^{-1}$ for $1\leq
i\leq n$.  Then $(b'_1,\ldots ,b'_n)\in Lg_n(A)$ (the proof is
similar to that of \cite{HeVa84}). We have
$$b_i-b'_i = [b_ib_{n+1}d-(c_i/K)b_{n+1}d_i][e/K+b_{n+1}d]^{-1}$$
Since $K\geq (1/\varepsilon') \Max(\| c_i\| : 1\leq i\leq n)$ and
$\| b_{n+1}d_i(e/K+b_{n+1}d')^{-1}\|\leq 1$ (this follows from 
$0\leq b_{n+1}d_i(e/K+b_{n+1}d)^{-1} \leq e$), we obtain
$$\| (c_i/K)b_{n+1}d_i(e/K+b_{n+1}d)^{-1}\|\leq\varepsilon' \;.$$ 
We also
have
$$0 \leq (b_ib_{n+1}dKv^{-1})^\ast(b_ib_{n+1}dKv^{-1})\leq
(b^{1/2}_0b_{n+1}dKv^{-1})^2\leq\varepsilon' \; .$$
Therefore $\| b_i-b'_i\| \leq \varepsilon'+\sqrt{\varepsilon'} \leq
2\sqrt{\varepsilon'} = \varepsilon$ and thus $tsr(A)\leq n$. The proofs of all
claims stated without proof are similar to the corresponding proofs in
\cite{HeVa84}. 
\myqed 
%
\section{Stable ranks and connectedness properties}

The characterizations of the Bass and the bilateral Bass stable ranks given above
relates them with a property of algebra or left module morphisms from $A$. We
will show that same property of algebra morphisms towards $A$
characterize the connected stable rank $csr(A)$. This will be obtained using some
Serre fibrations. A characterization for $acsr$ is also given. All connected
components in this section are in fact path-connected components.

\subsection{Serre fibrations}
By a result due to Michael \cite{Mic59}, an onto morphism $f : A\rightarrow B$ of
Banach algebras induces a Serre fibration $f^n : GL_n(A)\rightarrow GL_n(B)$.
 Also
\cite{CoSu87}, a (dense) onto morphism $f : A\rightarrow B$ of commutative unital
Banach algebras induces an (approximate) Serre fibration $f_n : Lg_n(A)\rightarrow
Lg_n(B)$. We will need corresponding facts for the more general situation of module
morphisms. 

Let $E,B$ and $X$ be Hausdorff topological spaces. A map $p : E\rightarrow B$ is
said to have the $HLP$--{\it{homotopy lifting property}\/}-- with
respect to $X$    if for every commutative diagram
$$
\begin{array}{ccc}
X & \stackrel{h}{\longrightarrow} & E\\
i\Big\downarrow && \Big\downarrow p\\
I\times X  & \sousfleche{\longrightarrow}{H} & B
\end{array}
$$
where $I = [0,1]$, $i(x) = (0,x)$, there exists a map
%
$F : I\times X\rightarrow E$ such that $Fi = h$ and $pF = H$. Then $p$ is called a
{\it{Serre fibration}\/} if $p$ has the $HLP$ with respect to every
cube $I^m$ ($m\geq 0$). We refer to \cite{Hu59} for more
information. The approximate $HLP$ and approximate Serre fibrations
were introduced in \cite{CoDu77} (they are called quasi-fibrations
by Corach and Su\'arez). We will not need this generalization in the
present paper.

Let $A$ be a unital Fr\'echet good algebra and let ${}_AM$ be a left Fr\'echet
module.  For
${\mathbf{m}} = (m_1,\ldots ,m_n)\in \Gen_n(M)$, let
$\Gen_n(M,{\mathbf{m}})$ be the connected component of $\Gen_n(M)$
containing ${\mathbf{m}}$. The following result was obtained by
Rieffel \cite[8.3]{Rie83} for projective left
$A$-modules $M$ and for Banach algebras $A$ (cf. also \cite{Tho91} for
the case of $C^\ast$-algebras). For Swan's problem, we will need to consider Fr\'echet
good algebras $A$ and the connected components of left $A$-modules of the form $M =
A/J$, for certain closed left ideals $J$. Not all modules of this type are
projective ; cf. \cite[pp. 346 ff.]{Hel93} for conditions implying
(equivalent to)  ``$A/J$ projective''.

\begin{Theo}\label{Theorem 3.1} 
Let $A$ be a unital Fr\'echet good algebra, let $M$ be
a left Fr\'echet module and $f : A\rightarrow M$ an onto module morphism. Then
$\Gen_n(M)$ is open in $M^n$ and we have $\Gen_n(M,{\mathbf{m}}) =
GL_n(A)_0{\mathbf{m}}$ for every ${\mathbf{m}}\in \Gen_n(M)$. 
\end{Theo}

{\bf Proof :} We show that $\Gen_n(M)$ is open in $M^n$. Let ${\mathbf{m}} =
(m_1,\ldots ,m_n) \in \Gen_n(M)$ be such that $b_1m_1 + \ldots +
b_nm_n = f(e)$, 
$b_i \in A$ ($i=1,\ldots ,n$). If $p_i = f(c_i)$, $c_i\in A$,
($i=1,\ldots ,n$) is sufficiently close of $m_i$, then
$f(\sum^n_{i=1}b_ic_i) = \sum^n_{i=1}b_if(c_i) =
\sum^n_{i=1}b_ip_i$ is close to $f(e) = \sum^n_{i=1}b_im_i$.
%
Since $A$ and $M$ are Fr\'echet spaces, the algebraic isomorphism $M\simeq A/\Ker
(f)$ is a topological isomorphism. This shows that there exists $j_0\in\Ker (f)$
such that $x+j_0$ is close to $e$, where $x = \sum^n_{i=1} b_ic_i$. This means
that $x+j_0$ is  invertible in $A : y(x+j_0) = e$ for a suitable $y\in A$. Thus
$$\sum^n_{i=1}yb_ip_i = yf(x) = yf(x+j_0) = f(e)\;,$$
yielding ${\mathbf{p}}\in \Gen_n(M)$.  

For the second claim, note that $GL_n(A)_0{\mathbf{m}}$ is contained in
$\Gen_n(M,{\mathbf{m}})$ since $GL_n(A)_0$ is path-connected. For the other
inclusion, we will show that $GL_n(A)_0{\mathbf{m}}$ is both open and closed
in $\Gen_n(M)$. Since $A$ is good, so is $M_n(A)$
\cite{Swa77}. Take an open neighborhood $U$ of $e$ in $GL_n(A)$.

The map $g : A^n\rightarrow M$ given by
$$g(a_1,\ldots ,a_n) = a_1m_1+\cdots + a_nm_n$$
is onto, since ${\mathbf{m}}\in\Gen_n(M)$. By the open mapping theorem
(for Fr\'echet spaces), $g$ is open. Therefore there exist an open
neighborhood $V'$ of $0$ in
$A^n$ such that ${\mathbf{m}} + V'\subset\Gen_n(M)$ and an open neighborhood
$V$ of $0$ in $A^n$ such that $g(V)\subseteq V'$ and $e + T\in U$ provided
that the rows of $T$ are in $V$. Then 
$$W = \{{\mathbf{v}} = (v_1,\ldots ,v_n)\in\Gen_n (M) :
{\mathbf{v}}-{\mathbf{m}}\in g(V)\}$$
is an open neighborhood of ${\mathbf{m}}$ in $\Gen_n(M)$. Let ${\mathbf{v}} =
(v_1,\ldots ,v_n)\in W$. For each $i=1,\ldots ,n$, there exists an
element 
$(b_{i1},b_{i2},\ldots ,b_{in})\in V$ such that
$\sum^n_{k=1}b_{ik}m_k = v_i-m_i$.
%
Then $(e+T){\mathbf{m}} = {\mathbf{v}}$, where $T = (b_{ij})\in M_n(A)$ and
$e+T\in U$. Thus $W$ is contained in $U{\mathbf{m}}$, showing, by
translation, that all orbits of $GL_n(A)$ in $\Gen_n(M)$ are open . Since
$GL_n(A)_0$ is open in  $GL_n(A)$, all orbits of $GL_n(A)_0$ in $\Gen_n(M)$
are open and thus also closed. Therefore $GL_n(A)_0{\mathbf{m}} =
\Gen_n(M,{\mathbf{m}})$.  \myqed

\smallskip

The following lemma is folklore for the case of $Lg_n(A)$. For $C^\ast$-algebras,
this follows easily from the characterization of $Lg_n(A)$ as elements
$(x_1,\ldots ,x_n) \in A^n$ with $\sum^n_{i=1}x^\ast_ix_i >
\varepsilon$ for a certain $\varepsilon > 0$. G. Corach
and F.D. Su\'arez \cite{CoSu92} recently proved that $Lg_n(C(X,A))$
is homeomorhic to $C(X,Lg_n(A))$ for compact
$X$ and Banach algebras $A$.

\begin{Lem}\label{Lemma 3.2} 
Let $X$ be a paracompact Hausdorff space and $A$ a
unital good topological algebra. Then $Lg_n(C(X,A)) \cong C(X,Lg_n(A))$. If $A$ is a
unital Fr\'echet good algebra and $f : A\rightarrow {}_AM$ is an onto left
$A$-module morphism, then $C(X,M)$ is a left module over $C(X,A)$ and
$\Gen_n(C(C,M))\cong C(X,\Gen_n(M))$. 
\end{Lem}

{\bf Proof :} Let $(f_1,\ldots ,f_n)\in Lg_n(C(X,A))$. Then
$$g_1f_1 +\ldots + g_nf_n = e\quad\hbox{\rm in}\quad C(X,A)$$
for a suitable $(g_1,\ldots ,g_n)\in C(X,A)^n$. This implies
$$(f_1,\ldots ,f_n) \in C(X,Lg_n(A)).$$

Let now $(f_1,\ldots ,f_n)\in C(X,Lg_n (A))$. Let $x\in X$ and
suppose that 
$a_1f_1(x) + \cdots + a_nf_n(x) = e$, for suitable $a_1,\ldots ,a_n$
in $A$. For $y$ close to $x$ in $X$, the element
$$z = a_1f_1(y) + \ldots + a_nf_n(y)$$
is close to $e$ in $A$ and thus invertible. Set $g_i(y) = z^{-1}a_i$ defined on a
neighborhood of $x$. Then $g_i$ are continuous (since $A$ is good) and 
$\sum_{i=1}^n g_i(y)f_i(y) = e$ on that neighborhood. Using a partition of unity
argument, we may glue together the local solutions $g_i$ to global ones, defined on
the paracompact space $X$. 

For the second part, note that $C(X,M)$ is a $C(X,A)$-module with the operation of
multiplication defined by $gh(t) = g(t)h(t), g \in C(X,A), h \in C(X,M), t \in
X$.  The
module morphism $f : A \rightarrow {}_AM$ induces a module morphism 
$f^X : C(X,A) \rightarrow C(X,M)$
by $f^X(g)(t) = f(g(t))$. 

We give here
only the argument concerning local identities. Let $(f_1, \ldots,
f_n) \in C(X,Gen_n(M))$ and suppose $a_1f_1(t) + \ldots + a_nf_n(t)
= f(e)$. According to the Bartle-Graves theorem for Fr\'echet spaces
\cite{Mic56}, there is a continuous  mapping
$\rho : M \rightarrow A$ such that $f(\rho(m)) = m,\; \;  (m \in M)$. Denoting
$c_i = \rho(f_i(t))$, $i = 1, \ldots, n$, we have $f(c_i) = f_i(t)$.
We obtain  
$$\sum_{i=1}^{n} a_if(c_i) = f(e)$$ 
and thus 
$$f(\sum_{i=1}^n a_ic_i) = f(e)\; .$$ 
Therefore  $\sum_{i=1}^n a_ic_i + j_0 =
e$, for a certain $j_0 \in Ker(f)$. For $s$ close to $t$, $b_i : = \rho(f_i(s))$ are
close to $c_i$, $i = 1, \ldots, n$, and 
$z : =  a_ib_i + j_0$ is close to e and thus
invertible. Then the functions $g_i(s) : = z^{-1}a_i$ are continuous on a
neighborhood of $t$ and, on that neighborhood,
$$f(e) = f(z^{-1}(\sum_{i=1}^n a_ib_i + j_0)) 
= \sum_{i=1}^n z^{-1}a_if_i(s)\; .$$
Thus  $\sum_{i=1}^n g_i(s)f_i(s) = f(e)$. A partition of unity arguments completes
the proof.  \myqed

\smallskip

\begin{Theo} \label{Theorem 3.3} 
Let $A$ be a Fr\'echet good algebra and let $f : A \rightarrow {}_AM$ be an onto
module morphism. Let ${\mathbf{m}} \in \Gen_n(M)$.Then 
$$T_{{\mathbf{m}}} : GL_n(A) \ni P \rightarrow P{\mathbf{m}} \in Gen_n(M)$$
is a Serre fibration. In particular, $T_{{\mathbf{m}}}$ is onto if and only
if the induced map $\pi_0(T_{{\mathbf{m}}}) : GL_n(A)/GL_n(A)_0 \rightarrow
\pi_0(\Gen_n(M))$ is onto. 
\end{Theo}

{\bf Proof :} To show that $T_{{\mathbf{m}}}$ has the $HLP$ with respect to 
$I^{0} = \{ 0 \}$, consider the commutative diagram
$$
\begin{array}{ccc}
\{ 0 \} & \stackrel{h}{\longrightarrow} & GL_n(A)\\
i\Big\downarrow && \Big\downarrow T_{{\mathbf{m}}}\\
I\times \{ 0 \}  & \sousfleche{\longrightarrow}{H} & \Gen_n(M)
\end{array}
$$
with $i(0) = (0,0)$. The onto morphism $f : A \rightarrow M$ induces an onto 
morphism, still denoted by $f$, between $C(I,A)$ and $C(I,M)$. By Lemma
\ref{Lemma 3.2},  we can think
of $H$ and $H(0,0)$ as elements of $\Gen_n(C(I,M))$; moreover, 
since $s \rightarrow H(st,0)$ is an arc in $\Gen_n(C(I,M)),$ 
they belong to the same
connected  component of
$\Gen_n(C(I,M))$. We have $H(0,0) = P{\mathbf{m}}$ for a certain 
$P = h(0) \in GL_n(C(I,A))$.
Here ${\mathbf{m}} \in C(I,M)$ is defined as the  constant
function equal to ${\mathbf{m}} \in M$ on $I$. 
By Theorem \ref{Theorem 3.1}, there exists a matrix 
$S \in GL_n(C(I,A))_0$ such that $H
= SP{\mathbf{m}}$. Note that $I$ and $I\times \{ 0 \}$ are
tacitly identified. Define $F \in C(I,GL_n(A))$ by 
$$F(t) : = S(t)P(t)(S(0)P(0))^{-1}h(0)\; .$$ 
Then $Fi = h$ and $T_{{\mathbf{m}}}F = H$, since $h(0){\mathbf{m}} = 
S(0)P(0){\mathbf{m}}$. Therefore $T_{{\mathbf{m}}}$ verifies $HLP$ with
respect  to $I^{0} = \{ 0\}$.

The commutative diagram
$$
\begin{array}{ccc}
X & \stackrel{h}{\longrightarrow} & GL_n(A)\\
i\Big\downarrow && \Big\downarrow T_{{\mathbf{m}}}\\
I\times X  & \sousfleche{\longrightarrow}{H} & \Gen_n(M)
\end{array}
$$
can be transformed, using Lemma \ref{Lemma 3.2}, into
$$
\begin{array}{ccc}
\{ 0 \} & \stackrel{g}{\longrightarrow} & GL_n(C(X,A))\\
i\Big\downarrow && \Big\downarrow T_{{\mathbf{m}}}\\
I\times \{ 0 \}  & \sousfleche{\longrightarrow}{G} & \Gen_n(C(X,M))
\end{array}
$$
where $g(0) = h$ and $G(t,0)(x) = H(t,x)$. It follows that
$T_{{\mathbf{m}}}$ has the $HLP$ with respect to all (paracompact) Hausdorff
spaces $X$. In particular, it is  a Serre fibration. \myqed

\smallskip

{\bf Remark :} The above proof is adapted from \cite{CoSu87}. It can be proved that
$T_{}$ is a principal fiber bundle  (cf. \cite{CoLa86b}). We refer to
\cite{CoLa86b}, \cite{Rie87},  \cite{Tho91}, \cite{Srd94} and \cite{Zha91} for
several instances of (Serre) fibrations.

\smallskip 

The following result was proved for morphisms of commutative Banch algebras by
Corach and Su\' arez \cite{CoSu87}.

\begin{Theo} \label{Theorem 3.4}
Let $A$ be a unital Fr\'echet good algebra. Every left onto
module  morphism $f : {}_AA \rightarrow {}_AM$ induces a Serre
fibration $f_n : Lg_n(A) \rightarrow \Gen_n(M)$. An onto algebra morphism $f :
A\rightarrow B$ between unital Fr\'echet good algebras induces a Serre fibration 
$f_n :
Lg_n(A) \rightarrow Lg_n(B)$.
\end{Theo}

{\bf Proof :} As above, it suffices to show that $f_n$ posses $HLP$ with respect
to  $I^{0} = \{ 0\}$. Consider the commutative diagram
$$
\begin{array}{ccc}
\{ 0\} & \stackrel{h}{\longrightarrow} & Lg_n(A)\\
i\Big\downarrow && \Big\downarrow f_{n}\\
I\times \{ 0\}  & \sousfleche{\longrightarrow}{H} & \Gen_n(M)
\end{array}
$$
and denote ${\mathbf{a}} : = h(0) \in Lg_n(A), {\mathbf{b}} : =
f_n({\mathbf{a}})
 = H(0,0) \in
Gen_n(M)$.  The mapping $T_{{\mathbf{b}}}$ is a Serre fibration. Therefore,
in the following commutative diagram
$$
\begin{array}{ccc}
\{ 0\} & \stackrel{g}{\longrightarrow} & GL_n(A)\\
i\Big\downarrow && \Big\downarrow T_{{\mathbf{b}}}\\
I\times \{ 0\}  & \sousfleche{\longrightarrow}{H} & \Gen_n(M)
\end{array}
$$
where $g(0) = e \in GL_n(A)$, there exists $F_1 : I\times \{ 0\} \rightarrow
GL_n(A)$ such that $F_1(0,0) = e$ and 
$T_{{\mathbf{b}}}(F_1(t,0)) = H(t,0), \; t \in I$. Define $F : I\times \{ 0\}
\rightarrow Lg_n(A)$ by $F(t,0) = T_{{\mathbf{a}}}(F_1(t,0))$. Then F is the
desired map. Indeed, $Fi = h$ and, using the equality $f_nT_{{\mathbf{a}}} = 
T_{f({\mathbf{a}})}$, we have $f_nF = H$.

The second part of the theorem follows from the first one, because every onto
algebra morphism $f : A \rightarrow B$ becomes an onto module morphism $f : {}_AA
\rightarrow {}_AB$, where ${}_AB$ is $B$ with the multiplication $a\ast b = f(a)b$.
Moreover, $\Gen_n({}_AB)$ is equal to $Lg_n(B)$. 	\myqed

\smallskip

This implies the following characterization for the (bilateral) Bass stable rank.

\begin{Cor} \label{CorSerre}
Let $A$ be a unital Fr\'echet good algebra. Then $Bsr(A) \leq
n$ if and  only if every left onto module morphism $f : {}_AA \rightarrow
{}_AM$ induces an onto mapping $\pi_0(f) : \pi_0(Lg_n(A))\rightarrow 
\pi_0(Gen_n(M))$. The condition $(bBsr)_n$ holds in $A$ if and
only if every onto algebra morphism $f : A \rightarrow B$
induces an onto mapping $\pi_0(f) : \pi_0(Lg_n(A)) \rightarrow \pi_0(Lg_n(B))$. 
\end{Cor} 

\subsection{Characterizations for csr} 
The following characterization of the
connected  stable rank is a generalization of a result due to Taylor \cite{Tay75}
and is a counterpart of the characterization of the (bilateral) Bass stable rank.
 
\begin{Theo} For a Fr\'echet good algebra $A$ the connected stable 
rank $csr(A)$ is  the least positive integer $m$ (or infinity if no
such number exists) such that for every  $n \geq m$ the following
property 
\begin{quote}
$(csr)_n$	\quad every onto unital algebra homomorphism $f : B \rightarrow A$ induces
an onto homomorphism  $f_n : Lg_n(B) \rightarrow Lg_n(A)$
\end{quote}
holds.
\end{Theo}

{\bf Proof :} We have to prove that the condition $(csr)_n$ is equivalent
to the connectedness of $Lg_n(A)$. 

Suppose first that $Lg_n(A)$ is connected. Let $f : B \rightarrow A$ be an onto
algebra morphism. It induces a map $\pi_0(Lg_n(B)) \rightarrow
\pi_0(Lg_n(A))$. Since $Lg_n(A)$ is connected, the map $\pi_0(f)$
is onto. Using Corollary \ref{CorSerre}, $f_n : Lg_n(B) \rightarrow Lg_n(A)$ is 
onto. Thus condition $(csr)_n$ holds. 

Suppose now that every onto morphism $B \rightarrow A$ induces an onto 
morphism $Lg_n(B) \rightarrow Lg_n(A)$. We want to prove that $Lg_n(A)$ is
connected. We take as $B$ the algebra of all continuous functions $f$ 
from $I = [0,1]$
into $A$ such that $f(0)$ is a multiple of the identity. This is a 
Fr\'echet algebra.
The map $f \rightarrow f(1)$ is an onto morphism from B onto A inducing a
morphism of $Lg_n(B)$ onto $Lg_n(A)$ (the surjectivity of $f$ is
immediate :  for $a
\in A$, take $f \in B$ as the linear map which joins $e$ with $a$). Now take 
${\mathbf{a}} \in
Lg_n(A)$ and consider ${\mathbf{f}} = (f_1,\ldots,f_n) \in Lg_n(B)$ such that
${\mathbf{a}} = {\mathbf{f}}(1)$. Then $\{ (f_1(t),\ldots,f_n(t)) : t \in
[0,1]\}$  is an
arc in $Lg_n(A)$ connecting ${\mathbf{a}}$ with an element of the form 
$(\lambda_1e, \ldots
,\lambda_ne) \in A^n$, with $(\lambda_1, \ldots ,\lambda_n) \in {\bf C}^n$. It 
follows that ${\mathbf{a}}$ can be connected with $e_n = (0, \ldots, 0,e)$,
for instance. Thus $Lg_n(A)$ is connected. \myqed

\smallskip

Another characterization of the connected stable rank for Banach algebras can be 
given in terms of the general stable rank. According to a result of Rieffel
\cite{Rie87}, the general stable rank of a Banach algebra $A$ is the least 
integer $k$
such that $Lg_n(A) = Lc_n(A)$ for all $n \geq k$, where $Lc_n(A)$ denotes the space
of last columns of invertible $n \times n$ matrices with entries in $A$ : 
$$Lc_n(A) = \{ Me_n : M \in GL_n(A)\} \; .$$ 
In Rieffel
\cite{Rie87} this result is stated for $C^{\ast}$-algebras, but the given 
proof works for all Banach algebras. We have
$gsr(A) \leq csr(A)$ (see \cite[p.328]{Rie83}). We define the
{\it {last columns connected stable rank}\/}  $Lccsr(A)$ as the least integer
$k$ such that for all $n \geq k$ the set $Lc_n(A)$ is connected. This is equivalent
to the fact that every element of $Lc_n(A),\;  n \geq k$, is the last column of a
matrix from $GL_n(A)_0$. The same proof of Theorem 2.1 in \cite{Rie87} 
works for all
Banach algebras, showing that $Lc_n(A)$ is connected if and only if the map from
$GL_{n-1}(A)$ to $GL_n(A)/GL_n(A)_0$ is surjective. 

\begin{Prop}
Let A be a unital Banach algebra. Then 
$$csr(A) = \max (gsr(A) ; Lccsr(A))\; .$$ 
\end{Prop}

{\bf Proof :} We have 
$$Lc_n(A)_0 = Lg_n(A)_0  = GL_n(A)_0e_n$$ 
(see \cite{Rie83} and \cite{Tho91} --
stated for $C^{\ast}$ - algebras --, or the proof of Theorem
\ref{Theorem 3.1}).  Let
$k \geq csr(A)$. Since $gsr(A) \leq csr(A)$, we have $k \geq gsr(A)$. For $n \geq
k$ we have $Lg_n(A) \supseteqÊ Lc_n(A) \supseteqÊ Lc_n(A)_0 = Lg_n(A)_0$. Since
$Lg_n(A)$ is connected, we have $Lg_n(A) = Lg_n(A)_0$ and thus $Lc_n(A) =
Lc_n(A)_0$ for all $n \geq k$. This implies $k \geq Lccsr(A)$. Therefore
$csr(A) \geq \max(gsr(A) ; Lccsr(A))$. Suppose now that $k \geq \max(gsr(A) ;
Lccsr(A))$. This means that for $n \geq k$ we have 
$Lg_n(A) = Lc_n(A)$ and $Lc_n(A)$ is
connected. Therefore $k \geq csr(A)$.	\myqed

\subsection{Characterizing acsr using tsr} 
Let ${\mathbf{B}}_r = \{ {\mathbf{x}} = (x_1, \ldots , x_r) \in {\bf R}^r :
x^2_1 + \cdots + x^2_r \leq 1\}$ be the closed ball in ${\bf R}^r$ and 
${\mathbf{S}}^{r-1} =
\partial {\mathbf{B}}_r$ the euclidean sphere. We will prove that $acsr(A) =
tsr(C(I,A))$ for all good normed algebras $A$. This generalizes the result due 
to V. Nistor
\cite{Nis86} mentioned in the first section.  
    
\begin{Lem}
Let  $A$  be a unital normed good algebra with 
jointly continuous multiplication, $n
\geq tsr(C({\mathbf{B}}_r,A))$  and  $f \in C({\mathbf{B}}_r,A^n)$ such that
$f(x) \in Lg_n(A)$  for all $x \in {\mathbf{S}}^{r-1}$. Then, for every
$\varepsilon > 0$, there exists  $g \in C({\mathbf{B}}_r ,Lg_n(A))$ such that 
$f\mid {\mathbf{S}}^{r-1}  = g\mid {\mathbf{S}}^
{r-1}$ and  $\|Êf - g \| < \varepsilon$. 
\end{Lem}
 
{\bf Proof :} We apply Theorrem \ref{Theorem 2.5}, $(i) \Rightarrow
(iii)$, with
$C({\mathbf{B}}_r,A)$ instead of $A$, with $C({\mathbf{S}}^{r-1},A)$ instead
of $B$ and with the canonical projection of $C({\mathbf{B}}_r,A)$ onto
$C({\mathbf{S}}^{r-1},A)$  instead of the
morphism $f$. Note that $C({\mathbf{B}}_r,A)$ is a good algebra if $A$
is.  Indeed, let $q
> 0$ be such that $\| a - e \| < q$ in $A$ implies that $a \in A^{-1}$. 
If $g \in
C({\mathbf{B}}_r,A)$ satisfies $\| g(t) - e \| < q$ for all $t \in
{\mathbf{B}}_r$, then $g(t) \in A^{-1}$. Then $g \in
C({\mathbf{B}}_r,A^{-1})$. It follows from the proof of Lemma \ref{Lemma 3.2}
that $g \in C({\mathbf{B}}_r,A)^{-1}$. The continuity of the inverse function
can be easily proved.  \myqed

\smallskip 

\begin{Theo} \label{Nis} 
Let $A$ be a good normed algebra. Let $n$ be a positive integer. Then $n$ belongs to
the absolute connected stable range of $A$ if and only if  $n
\geq tsr(C(I,A))$. 
Consequently,
$acsr(A) = tsr(C(I,A))$.
\end{Theo}

{\bf Proof :} The algebra $A$ is a quotient of $C(I,A)$ and thus $tsr(C(I,A))
\geq tsr(A)$  (see \cite[Theorem 4.3]{Rie83}). Let $n \geq tsr(C(I,A))$ and
let $V$ be an open, nonempty and connected subset of $A^n$. We will show that 
$V \cap
Lg_n(A)$ is connected. Let $x_0, x_1 \in V \cap Lg_n(A)$. Then there exists a
continuous function $f : I \rightarrow V$ such that $f(0) = x_0, f(1) = x_1$. 
Let $\varepsilon
> 0$ be less than the distance between the compact set $f(I)$ and 
the closed set $A^n
\setminus V$. Using the above Lemma for the interval $I$ we find a
continuous function $g : I \rightarrow Lg_n(A)$ such that $g(0) =
f(0) = x_0,  g(0) = f(0) = x_1$
and $\| g - f \| < \varepsilon$. Then $g$ takes its values in 
$\{ y : \dist (y,f(I)) < \varepsilon \}$
which is included, by hypothesis, in $V$. Thus $n$ is in the absolute connected
stable range.

For the converse implication, we can complete the proof as in \cite{Nis86}. 
\myqed

\section{Subalgebras with the same stable rank}
We prove in this section our results concerning Swan's problem. The main results,
including a proof of Theorem 1.1 cited in Introduction, are proved in section
\ref{main result}.

\subsection{Dense morphisms versus onto morphisms} 
Bass stable rank can be characterized using onto morphisms. We will
consider in this paragraph, as suggested by Corach and
Su\'arez \cite{CoSu87}, dense morphisms between topological algebras
instead of just onto morphisms. The same idea of looking at dense subrings of
topological rings is already known in $K$-theory. It goes back to the elementary
proof of the periodicity theorem of Atiyah and Bott and its very algebraic
exposition due to Bass. We refer to \cite{Bad94} for several related
references. 

The following result completes the characterization given in Theorem \ref{Theorem
2.2} for  the Bass stable rank and it is a compromise between the new ideology
(''dense morphisms``) and the old one (''onto morphisms``).

\begin{Theo} \label{Theorem 4.1} 
Let $A$ be a unital Fr\'echet good algebra. The following two
assertions are equivalent~: 
\begin{itemize}
	\item[(i)] The condition $(Bsr)_n$ holds, that is $Bsr(A) \leq n$ ; 
\item[(ii)] For every left $A$-module ${}_AM$, every onto module 
morphism $f : A \rightarrow  M$ induces a dense mapping $f_n :  Lg_n(A) \rightarrow
\Gen_n(M)$. 
\end{itemize}
\end{Theo}

{\bf Proof :} We have to prove that (ii) $\Rightarrow$ (i). Suppose that (ii)
holds and let $ {\mathbf{a}} = (a_1, \ldots , a_n , a_{n+1}) \in
Lg_{n+1}(A)$. Let $J$  be the closed
left ideal in $A$ generated by $a_{n+1}$. Let $M = A/J$ be the left normed
$A$-module and $\pi : A \rightarrow M$ the quotient  onto mapping. By (ii), the
mapping $\pi_n : Lg_n(A) \rightarrow \Gen_n(M)$ has a dense image. From the
definition of $M$ it follows that $\pi({\mathbf{a}})
= (\pi(a_1), \ldots, \pi(a_n)) \in \Gen_n(M)$. By Theorem 3.1, $\Gen_n(M)$ 
is open
in $M^n$ and thus $\pi({\mathbf{t}}) \in \Gen_n(M)$ if $\pi({\mathbf{t}})$ is
close to  $\pi({\mathbf{a}})$. 
 
By the density of $\pi_n$,
there exists ${\mathbf{x}} = (x_1, \ldots , x_n) \in Lg_n(A)$ such
that $\pi({\mathbf{x}})$ is very close to $\pi({\mathbf{a}})$. Consider the
commutative diagram
$$
\begin{array}{ccc}
\{ 0\} & \stackrel{h}{\longrightarrow} & Lg_n(A)\\
i\Big\downarrow && \Big\downarrow \pi_n\\
I\times \{ 0\}  & \sousfleche{\longrightarrow}{H} & \Gen_n(M)
\end{array}
$$
where $h(0) = {\mathbf{x}} \in Lg_n(A), i(0) = (0,0)$ and 
$$H(t,0) = \pi({\mathbf{a}}) + (1 -
t)(\pi({\mathbf{x}}) - \pi({\mathbf{a}}))$$ 
for $t \in I$. This shows that $H(t,0)$ is sufficiently close to 
$\pi({\mathbf{a}})$. Therefore $H(t,0) \in \Gen_n(M)$, for all 
$t \in I$. Since $\pi_n$ is a Serre fibration (Theorem 
\ref{Theorem 3.4}), there exists $F : I\times \{ 0\} \rightarrow Lg_n(A)$ such
that $Fi = h$ and 
$\pi_nF =
H$. Therefore 
$$\pi({\mathbf{a}}) = H(1,0) = \pi_n(F(1,0)) \in \pi_n(Lg_n(A))\; .$$ 
As in
the proof of Theorem \ref{Theorem 2.2}, it follows that $(a_1+c_1a_{n+1},
\ldots, a_n+c_na_{n+1}) \in Lg_n(A)$ for suitable $(c_1, \ldots, c_n) \in A^n$.
Thus $Bsr(A) \leq n$.  \myqed       

\smallskip

The above proof shows something more. Elements $(a_1, \ldots, a_n, a_{n+1}
) \in Lg_{n+1}(A)$ are said to be {\it{reducible}\/} (or stable \cite{Bas64}) 
when condition $(Bsr)_n$ holds for $(a_1, \ldots, a_n, a_{n+1})$, that is, when 
there exists $(c_1, \ldots, c_n) \in A^n$ such that $(a_1+c_1a_{n+1},
\ldots, a_n+c_na_{n+1}) \in Lg_n(A)$. The Bass stable rank of $A$ is no greater
than $n$ exactly when all $(n+1)$-tuples in $Lg_{n+1}(A)$ are reducible.

\begin{Theo} \label{Theorem 4.2}
	Let $A$ be a Fr\'echet good algebra. 
Consider $(a_1, \ldots, a_n, a_{n+1}) \in Lg_{n+1}(A)$ and let $J$ be the closed
left ideal in $A$ generated by $a_{n+1}$. The following two assertions are
equivalent : 
\begin{itemize}
	\item[(i)] The $(n+1)$-tuple $(a_1, \ldots, a_n, a_{n+1})$ is reducible ; 
\item[(ii)] The onto  module (quotient) morphism $\pi : A \rightarrow A/J$ 
induces a
dense mapping $\pi_n :  Lg_n(A) \rightarrow Gen_n(A/J)$.
\end{itemize} 
\end{Theo}

For the bilateral Bass stable rank we have the following characterization.
\begin{Theo} \label{Theorem 4.3} 
Let $A$ be a unital Fr\'echet good algebra. The following three 
assertions are equivalent : 
\begin{itemize}
	\item[(i)] The condition $(bBsr)_n$ holds  ; 
	\item[(ii)] For every  Fr\'echet good algebra $B$, every  continuous algebra 
morphism $f : A \rightarrow  B$, with $f(A)$ a dense and full subalgebra of 
$B$, induces a dense mapping $f_n :  Lg_n(A) \rightarrow Lg_n(B)$.
	\item[(iii)] Every onto  Fr\'echet good algebras morphism 
$f : A \rightarrow B$ induces 
a dense mapping $f_n :  Lg_n(A) \rightarrow Lg_n(B)$.
\end{itemize}
\end{Theo}	

{\bf Proof :} ''(i) $\Rightarrow$ (ii)`` Let 
${\mathbf{b}} = (b_1, \ldots , b_n) \in Lg_n(B)$ 
and let $d_i \in A, 1 \leq i \leq n$, with  $\sum_{i=1}^n d_ib_i = e$. 
Approximate sufficiently close $b_i$ and $d_i$ by $f(a_i)$ and $f(c_i)$, 
respectively, $a_i \in A, c_i \in A$, $i = 1, \ldots , n$. Then 
$\sum_{i=1}^n f(c_i)f(a_i)$ is close to $e = \sum_{i=1}^n d_ib_i$. Therefore 
$f(\sum_{i=1}^n c_ia_i)$ is invertible in $B$ and thus in $f(A)$. We 
obtain $f({\mathbf{a}}) = (f(a_1), \ldots , f(a_n)) \in Lg_n(f(A))$. Using 
Theorem \ref{Theorem 2.3}, the onto mapping $f : A \rightarrow f(A)$ induces 
an onto mapping $f_n : Lg_n(A) \rightarrow Lg_n(f(A))$. We can thus find 
${\mathbf{x}} = (x_1, \ldots , x_n) \in Lg_n(A)$ such that $f({\mathbf{x}}) = 
f({\mathbf{a}})$. Thus $f({\mathbf{x}})$ is close to ${\mathbf{b}}$ and 
${\mathbf{x}} \in Lg_n(A)$, yielding that $f_n : Lg_n(A) \rightarrow Lg_n(B)$ 
has dense range.

''(ii) $\Rightarrow$ (iii)`` Every onto algebra morphism $f : A
\rightarrow B$ is  a mapping such that $f(A) = B$ is a dense and
full subalgebra of $B$. $B$ is also  a Fr\'echet $Q$-algebra as it
can be easily proved.

''(iii) $\Rightarrow$ (i)`` is similar to the proof given in Theorem \ref{Theorem 
4.1}. \myqed

\smallskip

The {\it{dense stable rank}\/}  $dsr(A)$ of a Banach algebra 
$A$, as introduced by Corach and 
Su\'arez \cite{CoSu87}, is the least integer $n$ such that every
Banach  algebra dense 
morphism $f : A \rightarrow B$ induces a dense mapping 
$f_n : Lg_n(A) \rightarrow Lg_n(B)$. Theorem \ref{Theorem 4.3} 
implies the inequality $bBsr(A) \leq dsr(A)$. An affirmative answer for the 
problem ''dsr(A) = bBsr(A) (= Bsr(A)`` for all commutative $A$ 
would imply an affirmative answer to Swan's problem for commutative 
Banach algebras \cite{CoSu87}. It is known that 
$Bsr(H^{\infty}({\bf D})) = 1$
\cite{Tre92}, $dsr(H^{\infty}({\bf D})) = 1$ \cite{Sua94} and
$tsr(H^{\infty}({\bf D})) = 2$ \cite{Sua96}. 

\subsection{Reducible and bilateral reducible elements}
 
Consider a dense and full morphism $f : A \rightarrow B$ of two unital Fr\'echet 
good algebras $A$ and $B$. The inequality $Bsr(A) \leq Bsr(B)$ can be derived 
from the algebraic results of Swan \cite{Swa77}. We prove here a slightly 
more general result.

Recall that a morphism $f : A \rightarrow B$ is called \cite{CoSu87} $n$-{\it{
full}\/} if we have $f_n^{-1}(Lg_n(B)) \subseteq Lg_n(A)$ for the induced map 
$f_n : Lg_n(A) \to Lg_n(B)$. We call a morphism $f : A \rightarrow B$ 
{\it{near-unit full}\/} if there is an open neighborhood $V \subseteq B$ of 
$0$ such that $f(a) \in e + V$ implies $a \in A^{-1}$. If $f$ is full, 
then $f$ is also near-unit full since $B^{-1}$ is open. Note that Swan's more 
general question for topological rings was in fact raised for dense near-unit 
full morphisms.

\begin{Prop} \label{Proposition 4.4}
Let $A$ and $B$ be two unitaltopological algebras, $B$ being a $Q$-algebra, and 
$f : A \rightarrow B$ an algebra morphism with a dense image.
\begin{itemize}
\item[(i)] Suppose that $f$ is near-unit full. Then $f$ is $n$-full for 
every $n$.
\item[(ii)] Suppose that $Bsr(B) \leq n$ and $f$ is $n$-full. Then $Bsr(A) \leq 
n$.
\end{itemize}
\end{Prop} 

{\bf Proof :} (i) \quad Assume that $f$ is a dense near-unit full morphism. We will prove 
that $f$ is $n$-full for every $n$. Indeed, let 
$(f(a_1), \ldots , f(a_n)) \in Lg_n(B)$ with 
$$b_1f(a_1) + \cdots + b_nf(a_n) = e \; ; \; b_i \in B, i = 1, \ldots , n\; .$$
Let $V$ be that neighborhood of zero in $B$ which appears in the definition of 
near-unit fullness. Let $U \subseteq B$ be a neighborhood of $0$ 
such that $Uf(a_i) \subseteq (1/n)V$, for $i = 1, \ldots , n$. Approximate now 
$b_i$ by elements $f(c_i) \in f(A)$ such that $b_i - f(c_i) \in U$, $i = 1, 
\ldots , n$. Then
$$f(\sum_{i=1}^n c_ia_i) - e = \sum_{i=1}^n (f(c_i) - b_i)f(a_i) \in V$$
and thus $\sum_{i=1}^n c_ia_i$ is invertible in $A$. This implies that 
$f$ is $n$-full for all $n$.

(ii) \quad Suppose that $Bsr(B) \leq n$, $f$ is $n$-full, and let 
$(a_1, \ldots,a_{n+1}) \in Lg_{n+1}(A)$. This implies that $e \in A$ can be 
written as a linear combination by multiplying $a_i$ at the left by suitable 
elements. Applying $f$ to this equality we obtain $(f(a_1), \ldots
,f(a_{n+1})) 
\in Lg_{n+1}(B)$. Therefore, there exists ${\mathbf{d}} = (d_1,\ldots,d_n) 
\in B^n$ such that $f_n({\mathbf{a}}) + {\mathbf{d}}f(a_{n+1}) \in Lg_n(B)$. 
Approximate ${\mathbf{d}}$ by $f({\mathbf{c}}), {\mathbf{c}} \in A^n$. Then 
$f({\mathbf{a}}) + f({\mathbf{c}})f(a_{n+1})$ is close to $f({\mathbf{a}}) + 
{\mathbf{d}}f(a_{n+1}) \in Lg_n(B)$. Since $Lg_n(B)$ is open in $B^n$, 
${\mathbf{c}}$ can be chosen such that $f({\mathbf{a}}) +
f({\mathbf{c}})f(a_{n+1})  \in Lg_n(B)$. Since $f$ is $n$-full,
${\mathbf{a}} + {\mathbf{c}}a_{n+1}  \in Lg_n(A)$. Therefore $Bsr(A) \leq n$.
\myqed

\smallskip

For the reverse inequality ''$Bsr(B) \leq Bsr(A)$`` we can prove the following
partial result.  We will use in the proof a modification of 
$\Gen_n(M)$. The {\it{dense}\/} (or topological) 
$n$-{\it{generator space}\/} of
the left $A$-module $M$, denoted by $\Gen^{\ast}_n(M)$, is the set of elements 
$(v_1, \ldots,v_n) \in M^n$ such that
$Av_1 + \ldots + Av_n$  is dense in $M$. We have $\Gen_n(M) \subseteq 
\Gen^{\ast}_n(M)$. For $M = {}_AA$, 
$\Gen_n(A) = \Gen^{\ast}_n(A) = Lg_n(A)$.

\begin{Theo} \label{Theorem 4.5} 
Let $A$ and $B$ be two unital Fr\'echet 
good algebras, $Bsr(A) \leq n$, and $f : A \rightarrow B$ a dense
algebra morphism such that $f(A)$ is a full subalgebra of $B$. 
Let $(b_1, \ldots, b_{n+1}) \in  Lg_{n+1}(B)$ and let $J$ be 
the closed left ideal in $B$ generated by $b_{n+1}$. If $J_0
= J \cap f(A)$ is dense in $J$, then 
$(b_1, \ldots, b_{n+1})$ is reducible in $B$.  
\end{Theo}

{\bf Proof :} Note that if $n = \infty$ there is nothing 
to prove. Assume $n < \infty$.

\smallskip

\noindent {\sc Step 1} ({\bf the module} ${}_AM$)

\smallskip

Denote $M = B/J$ the left module over $B$ and consider the
quotient projection 
$\pi : B
\rightarrow {}_BM$ and the map $A \times M
\rightarrow M$ denoted by $\ast$ and defined by 
$$a\ast v = f(a)v, \; \; a \in
A, v \in M \; .$$ 
Then
${}_AM$ with its internal additive operation 
(inherited from ${}_BM$ as a
$B$-module) and the external operation $\ast$ becomes a normed left
$A$-module. Indeed, we have 
$$a\ast (v_1+v_2) = f(a)(v_1+v_2) = f(a)v_1 + f(a)v_2 = a\ast v_1 +
a\ast v_2
\; ;$$
the equality $(a_1+a_2)\ast v  = a_1\ast v+ a_2\ast v$ can be proved
in a similar fashion. Since $f$ is an algebra morphism, we have
$$(a_1a_2)\ast v = f(a_1a_2)v = f(a_1)(f(a_2)v) = a_1\ast (f(a_2)v) 
= a_1\ast (a_2\ast v) \; .$$
Consider $\pi f : A
\rightarrow {}_AM$. This mapping is a morphism of left
modules. The linearity is clear and also
$$(\pi f)(ab) = \pi(f(a)f(b)) = f(a)\pi(f(b)) = 
a\ast (\pi f)(b)\; .$$
The map $\pi f$ has a dense range since $f$ has a dense range 
and $\pi$ is onto. 

\smallskip

\noindent {\sc Step 2} ({\bf the map} $(\pi f)_n$ 
{\bf is dense})

\smallskip

We show that the map $\pi f : A \rightarrow {}_AM$ 
induces a dense mapping $(\pi f)_n : Lg_n(A) \rightarrow
\Gen^{\ast}_n({}_AM)$. Note the use here of the dense generator
space.

Let $v = (v_1, \ldots , v_n) \in
\Gen^{\ast}_n({}_AM)$. Then there exist $y_i \in A, 1 \leq i 
\leq n$, such that $\sum_{i=1}^n y_iv_i$ is very close to 
$\pi f(e)$. Since the range of $\pi f$ is dense in $M$, 
we can well approximate $v_i$ by elements 
$\pi f(x_i)$, $x_i \in A$, $i = 1, \ldots , n$, such that 
$\sum_{i=1}^n y_i \ast (\pi f)(x_i)$ is close to $\pi f(e)$. 

Denote $a = \sum_{i=1}^n y_ix_i$. Since $\pi f(a)$ is close to 
$\pi f(e)$, there exists $j \in J$ such that $f(a) - f(e) + j$ 
is small. But
$J_0$ is dense in $J$~; thus there exists $j_0 = f(t) \in 
J \cap f(A), \; t \in A$, such that $f(a) - e + j_0$ is 
sufficiently small. This implies that $f(a) + 
f(t)$ is invertible in $B$ and thus in the full subalgebra $f(A)$.
Hence, there is $c \in A$
such that $f(c)(f(a) + f(t)) = e$. Consider this equality in $B$
and apply the module morphism $\pi$. We obtain 
$f(c)\pi (f(a)) =
\pi(e)$, yielding
$\pi (f(a)) \in \Gen_1(X)$, where ${}_AX$ is the 
submodule $(\pi f)(A)$
of ${}_AM$. 

We have 
$$\sum_{i=1}^n (cy_i)\ast (\pi(f(x_i)) = 
c\ast(\pi\circ f)(\sum_{i=1}^n y_ix_i) =
f(c)(\pi f)(\sum_{i=1}^n y_ix_i) =
\pi(e)\; .$$
This shows that $(\pi\circ f)({\mathbf{x}}) \in \Gen_n({}_AX)$, 
${\mathbf{x}} = (x_1, \ldots, x_n) \in A^n$. Since
$Bsr(A) \leq n$, the onto module morphism $\pi f : A \rightarrow
{}_AX$ induces, by Theorem 2.2, an onto mapping $(\pi f)_n :
Lg_n(A) \rightarrow \Gen_n({}_AX)$. Therefore there exists
${\mathbf{b}} = (b_1,\ldots, b_n) \in Lg_n(A)$ such that $(\pi \circ
f)({\mathbf{b}}) = (\pi f)({\mathbf{x}})$. Hence $(\pi \circ
f)({\mathbf{b}})$  is close to ${\mathbf{v}}$ and the
mapping $(\pi f)_n : Lg_n(A) \rightarrow 
\Gen^{\ast}_n({}_AM)$ has
a dense range. 
\smallskip

\noindent {\sc Step 3} ({\bf conclusion})

\smallskip

We have 
\begin{eqnarray*}
\Gen^{\ast}_n({}_AM) &=& \{(v_1, \ldots, v_n) \in M^n :  A\ast v_1 +
\cdots +  A\ast v_n \mbox{ dense in } M\}\\
&=& \{(v_1, \ldots, v_n) \in M^n :  f(A)v_1 + \cdots + f(A)v_n 
\mbox{ dense in } M\}\\
&= & \{ (v_1, \ldots, v_n) \in M^n :  Bv_1 + \cdots + Bv_n
\mbox{ dense in } M \}\\ 
& &(\mbox{ since }f(A) \mbox{ is dense in } B )\\
&= & \Gen^{\ast}_n({}_BM) \; .
\end{eqnarray*}
We also have 
$$(\pi f)(Lg_n(A)) \subseteq  \pi(Lg_n(B)) \subseteq \Gen_n({}_BM) 
\subseteq \Gen^{\ast}_n({}_BM) =  \Gen^{\ast}_n({}_AM)\; .$$

Indeed, the first inclusion follows because $f(Lg_n(A)) \subseteq
Lg_n(B)$.  For the second one we use the fact that $\pi$ is an onto
morphism of left modules between $B$ and ${}_BM$, the last inclusion
is clear,  while
the given equality was proved before. Using Step 2, we get the
density of
$\pi(Lg_n(B))$  in
$\Gen_n(B/J)$. By Theorem \ref{Theorem 4.2}, $(b_1,
\ldots, b_{n+1}) \in Lg_{n+1}(B)$ is reducible. 	\myqed

\smallskip

\begin{Cor} \label{Cor 4.6}
Let $A$ and $B$ be two unital Fr\'echet good algebras, $Bsr(A) \leq 
n$, and $f : A \rightarrow B$ a dense algebra morphism such that 
$f(A)$ is a full subalgebra of $B$. Then all elements $(b_1,
\ldots, 
b_{n+1}) \in Lg_{n+1}(B)$ with $b_{n+1} \in f(A)$ are reducible 
in $B$.
\end{Cor}

{\bf Proof :} Let $b_{n+1} = f(c)$, $c \in A$. The closed left ideal 
$J$ generated by $b_{n+1}$ is the closure in $B$ of $Bf(c)$. 
Then 
$$f(A)f(c) \subseteq J_0 \subseteq J = \overline{f(A)f(c)}$$
since $f(A)$ is dense in $B$. Here $J_0 = J \cap f(A)$. It follows 
that $\overline{J_0} = J$. \myqed
\smallskip

An element $(b_1,
\ldots, b_{n+1}) \in Lg_{n+1}(B)$ is said to be {\it{bilateral 
reducible}\/} (in $B$) whenever there exist $(c_1,
\ldots, 
c_{n})$, $(d_1,\ldots, d_{n})$ in $B^n$ such that 
$$(b_1 + c_1b_{n+1}d_1,
\ldots, 
b_n + c_nb_{n+1}d_n) \in Lg_n(B)\; .$$
Theorem \ref{Theorem 4.2} has an analogous counterpart for 
bilateral reducible elements.

\begin{Theo} \label{Theorem 4.7}
Let $f : A \rightarrow B$ be a dense and full morphism of 
unital Fr\'echet good algebras.
\begin{itemize}
\item[(i)] We have $bBsr(A) \leq bBsr(B)$ ;
\item[(ii)] Suppose $bBsr(A) = n$ and let 
$(b_1,\ldots, b_{n+1}) \in Lg_{n+1})(B)$. Let $J$ be the closed
two-sided ideal in $B$ generated by $b_{n+1}$. If 
$J_0 = J \cap f(A)$ is dense in $J$, then 
$(b_1,\ldots, b_{n+1})$ is bilateral reducible in $B$. In
particular, all elements $(b_1,
\ldots, b_{n+1}) \in Lg_{n+1})(B)$ with $b_{n+1} \in f(A)$ are 
bilateral reducible in $B$.
\end{itemize}
\end{Theo}

{\bf Proof :} (i) \quad Suppose condition $(bBsr)_n$ holds in 
$B$ and let $(a_1,\ldots , a_{n+1}) \in Lg_{n+1}(A)$. Then 
$f(a_1),\ldots,f(a_{n+1})) \in Lg_{n+1}(B)$. Therefore 
$f({\mathbf{a}}) + {\mathbf{c}}f(a_{n+1}){\mathbf{d}} \in Lg_n(B)$ 
for suitable ${\mathbf{c}} = (c_1,\ldots,c_n)$, ${\mathbf{d}} =
(d_1,\ldots,d_n) \in B^n$. By the density of $f(A)$ in $B$, we 
can approximate ${\mathbf{c}}$ and ${\mathbf{d}}$ by $f({\mathbf{x}})$ 
and $f({\mathbf{y}})$, respectively, ${\mathbf{x}},{\mathbf{y}} 
\in A^n$, such that $f({\mathbf{a}}) +
f({\mathbf{x}})f(a_{n+1})f({\mathbf{y}})$ is close to $f({\mathbf{a}}) + 
{\mathbf{c}}f(a_{n+1}){\mathbf{d}} \in Lg_n(B)$. Since $Lg_n(B)$ is 
open in $B^n$, we obtain $f({\mathbf{a}}) +
f({\mathbf{x}})f(a_{n+1})f({\mathbf{y}}) \in Lg_n(B)$. Then 
$$\sum_{i=1}^n w_i(f(a_i + f(x_i)f(a_{n+1})f(y_i) = e$$
for suitable ${\mathbf{w}} = (w_1,\ldots,w_n) \in B^n$. Let 
${\mathbf{v}} = (v_1,\ldots,v_n) \in A^n$ be such that
$f({\mathbf{v}})$ is sufficiently close to ${\mathbf{w}}$ in the
$Q$-algebra $B$. Then $f(\sum_{i=1}^n v_i(a_i + x_ia_{n+1}y_i))$ is 
invertible in $B$. Since $f$ is full, $\sum_{i=1}^n v_i(a_i +
x_ia_{n+1}y_i)$ is invertible in $A$. This implies that condition
$(bBsr)_n$ holds in $A$. 

The proof of (ii) is similar to that of Theorem \ref{Theorem 4.5}. 
Note also that here we can use $\Gen_n(M)$ instead of 
$\Gen_n^{\ast}(M)$. \myqed

\subsection{Swan's problem for subalgebras of $C^{\ast}$-algebras} 
\label{main result}

The two last results can be succesfully applied for subalgebras of 
$C^{\ast}$-algebras for which the condition $\overline{J_0} = J$ in Theorem
\ref{Theorem 4.5} holds. 

Consider the following class of $\ast$-subalgebras of C$^{\ast}$-algebras.

{\bf Definition : } (\cite{KiSh93}) Let $C$ be a $\ast$-algebra wich is a 
dense subalgebra of the $\ast$-algebra $D$ with common identity $e$. If $D$ is 
commutative, $C$ is called {\it normal}\/ if the algebras of functions 
$\{ x(s) : x \in C \}$ is normal on the space $S$ of maximal ideals of $D$, 
that is there exists $x \in C$ such that $x(s) = 0$ on $T_1$ and $x(s) = 1$ 
on $T_2$ for every disjoint closed subsets $T_1$ and $T_2$ in $S$. In the
general case, $C$ is called {\it locally normal}\/ if for every $x \in C$,
there exists a commutative Banach $\ast$-subalgebra $D(x)$ in $D$ which
contains 
$e$ and $x$ and such that $C(x) = D(x) \cap C$ is a dense normal subalgebra of 
$D(x)$.

\smallskip

The class of dense locally normal $\ast$-subalgebras of $C^{\ast}$-algebras
contains the class of dense $\ast$-subalgebras of $C^{\ast}$-algebras
closed under $C^{\infty}$-functional calculus of self-adjoint elements
\cite{KiSh93}. The Wiener algebra shows that the inclusion is strict \cite{KiSh93}.

It was shown in the same paper \cite[Th. 13]{KiSh93} that, assuming $C$ locally
normal in the  $C^{\ast}$-algebra $D$, $J \cap C$ is dense in $J$ for all two-sided
ideals in $D$. This and theorem \ref{Theorem 4.5} imply
the following result : 

\begin{Cor} \label{Cor 4.8}
Let $f : A \rightarrow B$ be a dense and full $\ast$-morphism between the
unital Fr\'echet good algebra $A$ and a $C^{\ast}$-algebra $B$. Suppose $f(A)$
is locally normal in $B$. Then  $bBsr(A) = bBsr(B)$.
\end{Cor}

The same result holds for the Bass stable rank.

\begin{Theo} \label{Theorem 4.9}
Let $f : A \rightarrow B$ be a dense and full $\ast$-morphism between the
unital Fr\'echet $Q$-algebra $A$ and a $C^{\ast}$-algebra $B$. Suppose $f(A)$
is locally normal in $B$. Then 
$Bsr(A) = Bsr(B)$.
\end{Theo}

{\bf Proof :} By Proposition \ref{Proposition 4.4} $Bsr(A) \leq Bsr(B)$. For 
the reverse inequality,  suppose that 
$Bsr(A) = n < \infty$. According to the result of Herman and Vaserstein
\cite{HeVa84}, one has $Bsr(B) = tsr(B)$. We show that $tsr(B) \leq n$.

Let $\varepsilon$ be a
positive number,  $0 < \varepsilon < 4$, and set $\varepsilon' =
\varepsilon^2/16 < 1$. Let $(b_1, \ldots , b_n) \in B^n$ be given. 
By the density
of $f$, there exists $a_1, \ldots, a_n \in A$ such that $\| b_i - f(a_i) \| \leq
\varepsilon /2$, $i = 1, \ldots, n$. Set $b_0 =  \sum_{i=1}^n
f(a_i)^{\ast}f(a_i) \in f(A)$ and $x = e - b_0/\varepsilon'$. Note that the
self-adjoint element $x$ belongs to the commutative $C^{\ast}$- algebra $C$
generated by $e$ and $b_0$. Let  $B(x)$ be a commutative $C^{\ast}$-subalgebra
containing $e$ and $x$ such that  $f(A) \cap B(x)$ is normal. Let $S$ be the
maximal ideal space of $B(x)$ and  $$ T_1 = \{ s \in S : x(s) \leq \frac{1}{3} \}
\; \; ; \; \;  T_2 = \{ s \in S : x(s) \geq \frac{2}{3} \}\; .$$
There exists $y \in f(A) \cap B(x)$ such that $y(s) = 0$ for $s \in T_1$ and 
$y(s) = 1$ for $s \in T_2$. Replacing eventually $y$ by $y^{\ast}y$, we 
can assume that $y$ is positive. If $b_0$ vanish 
in a point of the maximal ideal of $C$, then $x$ in this point equal 1 and 
thus $y$ does not vanish there. This shows that $b_0 + y$ is invertible in $B$ 
and thus $(f(a_1), \ldots , f(a_n), y) \in Lg_{n+1}(B)$.
Since $Bsr(A) = n$ and $y \in f(A)$, this $(n+1)$-tuple is, by Corollary
\ref{Cor 4.6}, reducible. Therefore, there exists $(c_1, \ldots ,c_n) 
\in B^n$
such that $(f(a_1) + c_1y, \ldots, f(a_n) + c_ny) \in Lg_n(B)$. Let  $K \geq
(1/\varepsilon')\max (\| c_i \|: 1 \leq i \leq n)$ be fixed and define $v = e +
Ky$. Then $v \geq e$, so it is invertible. Let  $b'_i =  (f(a_i) + c_iy)v^{-1}$
for $1 \leq i \leq n$. Then, similarly as in \cite{HeVa84} and in
Theorem~\ref{Theorem 2.10}, one has $(b'_1 , \ldots , b'_n) \in Lg_n(A)$ and 
$$\| b_i - b'_i \| \leq \| b_i - f(a_i) \| + \| f(a_i) - b'_i \| \leq 
\varepsilon /2 + \varepsilon' + \sqrt{\varepsilon'} \leq  \varepsilon \; ,$$ 
for $i =
1, \ldots, n$. Thus $Bsr(B) = tsr(B) \leq n$. \myqed

\smallskip

This proves in particular Theorem 1.1. 

Several
examples of locally normal subalgebras or of subalgebras closed under
$C^{\infty}$-functional calculus of self-adjoint elements can be found in 
\cite{Co96}, \cite{BlCu91}, \cite{KiSh93}, \cite{KiSh94}. We state here only two
consequences of the above results.

Consider a closed derivation $\delta$ from a dense subalgebra $D(\delta) =
\mbox{ Domain }(\delta)$ of the unital $C^{\ast}$-algebra $A$ into $A$. 
For every $n \geq 2$, set
$$D(\delta^n) = \{x \in D(\delta) : \delta^k(x) \in D(\delta) \mbox{ for } 
1 \leq k \leq n-1 \}$$
and $D^{\infty}(\delta) =   \cap_{n=1}^{\infty} D(\delta^n)$.

\begin{Cor} \label{der}
Let $A$ be a unital unital $C^{\ast}$-algebra and $\delta$ a closed derivation. 
Then $Bsr(D(\delta)) = Bsr(A)$. If $D(\delta^n)$ is dense in $A$ for a certain 
$1 < n \leq \infty$, then $Bsr(D(\delta^n)) = Bsr(A)$. 
\end{Cor}

{\bf Proof :} This follows from Theorem \ref{Theorem 4.9} and the fact that the
subalgebras $D(\delta)$ and $D(\delta^n)$ of $A$ are closed under 
the $C^{\infty}$-functional calculus of self-adjoint elements and thus locally
normal (cf. for instance \cite{KiSh93}). \myqed

\smallskip 

We use now Theorem \ref{Theorem 4.9} to prove the equality of the Bass 
stable rank of a subalgebra of the reduced group $C^{\ast}$-algebra
$C^{\ast}_r(G)$ of a discrete group $G$, which is rapidly decaying \cite{CoMo90}.

Recall the following definition. A length function  $L$  on a topological group 
$G$ is a continuous nonnegative real valued function on $G$, such that $L(1) = 0$,
$L(g^{-1}) = L(g)$ and $L(gh) \leq L(g) + L(h)$, for all $g$ and $h$ in $G$. If $G$
is a discrete group with a given length function $L$, we define 
$$H_L^{\infty}(G) : = \{ \varphi \in \ell^2(G) :  \sum_{g\in G} |\varphi(g)|
{}^2(1 + L(g))^{2k} <  \infty ,  \mbox{ for all } k > 0 \}\; .$$ %
The group $G$ is said \cite{CoMo90}, \cite{Jol90} to be {\it rapidly decaying\/}  
if there exists a length function $L$ such that $H_L^{\infty}(G)$ is contained in
$C^{\ast}_r(G)$. $H_L^{\infty}(G)$ is then a Fr\'echet algebra with respect to the 
seminorms 
$$\| \varphi \|_k = ( \sum_{g\in G} |
\varphi(g)| {}^2(1 + L(g))^{2k})^{1/2}\; , \quad k \geq 1, \quad \varphi 
\in H_L^{\infty}(G). $$ 
The algebra $H_L^{\infty}(G)$ plays the role of smooth ($C^{\infty}$) functions : 
if $G = {\bf Z}^k$,  then $C^{\ast}_r(G) = C({\bf T}^k)$ and $f \in C({\bf T}^k)$ is
smooth if and only if its Fourier transform (a function on ${\bf Z}^k$) is of rapid
decay (cf. \cite{Jol90}). 

\begin{Cor} Let $G$ be a discrete rapidly decaying group with group length $L$. 
Then $Bsr(H_L^{\infty}(G)) = tsr(H_L^{\infty}(G)) = Bsr(C^{\ast}_r(G)) =
tsr(C^{\ast}_r(G))$. For $tsr$ we consider $H_L^{\infty}(G)$ endowed with the
topology endowed from $C^{\ast}_r(G)$. 
\end{Cor}

{\bf Proof :} This fits into the framework of Corollary \ref{der} via a contruction
due to Connes and Moscovici \cite{CoMo90} and a result from \cite{Ji92}. Indeed, the
operator $D : \ell^2(G) \rightarrow \ell^2(G)$ given by $(D\xi)(g) = L(g)\xi(g)$,
$g \in G$, is a self-adjoint and closed unbounded linear operator. Then $\delta a =
i[D,a]$ defines a closed, unbounded derivation from $B(\ell^2(G))$ into 
$B(\ell^2(G))$ (cf. \cite{CoMo90}, \cite{Ji92}). Consider the restriction of
$\delta$ to $C^{\ast}_r(G)$, denoted also by $\delta$. Then, according to Theorem
1.3 in \cite{Ji92}, $H_L^{\infty}(G) =   \cap_{n=1}^{\infty}
\mbox{Domain}(\delta^n)$ which is dense in $C^{\ast}_r(G)$. Therefore we can apply
Corollary \ref{der}. For the topological stable rank one equality follows from
\cite{HeVa84}. The other equality can be proved directly or using the results of
the next section.   \myqed
 
\subsection{Swan's problem for $tsr$ and $acsr$} 
The following simple
result shows that $tsr(B) \leq tsr(A)$ holds for all dense
morphisms $f : A \rightarrow B$.

\begin{Prop} 
Let $A$ and $B$ be two unital topological algebras and $f : 
A \rightarrow B$ a continuous algebra morphism with a dense image.
Then $tsr(B) \leq tsr(A)$. 
\end{Prop}

{\bf Proof :} Let $tsr(A) = n < \infty$, ${\mathbf{b}} = 
(b_1,\ldots,b_n) \in B^n$
and let  $V$ be an open neighborhood of $0$ in $B^n$. Since $f$ is
continuous, there is an open neighborhood $U$ of $0$ in $A^n$ such
that $f(U) \subseteq (1/2)V$. By the density of $f$, there
exists ${\mathbf{a}} = (a_1,
\ldots,a_n) \in A^n$ such that ${\mathbf{b}} - f({\mathbf{a}})$ belongs
to $(1/2)V$. Because $tsr(A) = n$, there exists 
${\mathbf{x}} = (x_1, \ldots, x_n)
\in Lg_n(A)$ such that ${\mathbf{a}} - {\mathbf{x}} \in U$. We obtain
$${\mathbf{b}} - f({\mathbf{x}}) = {\mathbf{b}} - f({\mathbf{a}}) +
f({\mathbf{a}}) - f({\mathbf{x}}) \in  (1/2)V + f(U) \subseteq V\;
.$$    
Therefore $tsr(B) \leq
tsr(A)$, since $f({\mathbf{x}}) \in Lg_n(B)$ (apply $f$ to 
an equality of the form $e = \sum_{i=1}^n y_ix_i$). \myqed

\smallskip

The reverse inequality is also true for some topological
algebras $A$ and $B$, provided $f$ satisfies an  additional
condition, stronger than ''$f$ is near-unit full``. We  call a 
continuous morphism $f : A \rightarrow B$ of unital topological algebras
{\it{strong near-unit full}\/} if for every open neighborhood $U$ of
$0$ in $A$, there exists an open neighborhood $V$ of $0$ in $B$ such
that $e - f(a) \in V$ implies $a$ invertible in $A$ with an inverse
$a^{-1}$ such that $e - a^{-1} \in U$.

A strong near-unit full
morphism is near-unit full. A dense strong near-unit full morphism
is full. It is easy to see that if $A$ is a dense and full subalgebra
of the good algebra $B$, then the full inclusion morphism $i : A
\rightarrow B$ is actually strong near-unit full by the continuity of
the inverse function. If $A \subseteq B$ and the inclusion is strong
near-unit full, then $A$ has a topology equivalent to $B$'s topology
as long as the inversion is continuous in $A$. This remark was 
communicated to the author by
L.B. Schweitzer. Indeed,
if $(a_n)$ is a sequence in $A$ converging to $0$ in topology of
$B$, then $e - a_n$ is eventually invertible in $B$ and so in $A$ (by
the fullness). Also $(e - a_n)^{-1}$ tends to $e$ in $A$, so $e -
a_n$ tends to $e$. Thus $a_n$ tends to $0$ in $A$.

For simplicity, the following result is stated for good normed algebras.

\begin{Theo} \label{snuf}
Let $A$ and $B$ be two unital good normed algebras and $f : 
A \rightarrow B$ a strong near-unit full algebra morphism with
dense image. Then $tsr(A) = tsr(B)$. In particular, this holds if
$A$ is a dense and full subalgebra of the good algebra $B$,
endowed with the topology of $B$.
\end{Theo}

{\bf Proof :}  Let $tsr(B) = n$ and let ${\mathbf{a}} = (a_1, \ldots,
a_{n+1})
\in Lg_{n+1}(A)$. Let $\varepsilon > 0$. We will show, using the
characterization of $tsr$ given in Theorem \ref{Theorem 2.6},
that $tsr(A) \leq n$. Set $\varepsilon' =
\varepsilon/(1+\varepsilon) < 1$. Then $(f(a_1),
\ldots, f(a_{n+1})) \in Lg_{n+1}(B)$ and, by Theorem
\ref{Theorem 2.6}, there exists ${\mathbf{d}} = (d_1,
\ldots, d_n) \in B^n$ such that $f({\mathbf{a}}) +
{\mathbf{d}}f(a_{n+1}) \in Lg_n(B)$ and $\| d_if(a_{n+1}) \| <
\delta(\varepsilon')/2$, $i = 1, \ldots, n$. Here
$\delta(\varepsilon')$ is a positive number such that $\| e -
f(a)\| 
< \delta(\varepsilon')$ in $B$ implies $a \in A^{-1}$ and $\| e -
a^{-1} \| < \varepsilon'$. Thus we can find ${\mathbf{y}} = (y_1,
\ldots , y_n)
\in B^n$ such that $\sum_{i=1}^n y_i(f(a_i)+d_if(a_{n+1})) = e$.

Suppose first that $f(a_{n+1}) \neq 0$. Let 
$M = \max\{\| f(a_i) \| : i
= 1,
\ldots , n\}$, $N = \max\{
\| y_i \|: i = 1,\ldots,n\}$ and set 
$$\varepsilon'' = \min\{1 ,
\frac{\delta(\varepsilon')}{2\|f(a_{n+1})\|} , 
\frac{\delta(\varepsilon')}{2n(N+1)\|f(a_{n+1})\|} ,
\frac{\delta(\varepsilon')}{n(2M+\delta(\varepsilon')}\}\; .$$ 
Choose ${\mathbf{c}} = (c_1, \ldots , c_n)$, ${\mathbf{x}} =
(x_1,\ldots, x_n)
\in A^n$ such that $\| d_i - f(c_i) \| < \varepsilon''$, $\| y_i -
f(x_i)
\| < \varepsilon''$. Then
\begin{eqnarray*}	
\| f(\sum_{i=1}^n x_i(a_i + c_ia_{n+1})) - e \|  & 
	\leq & \sum_{i=1}^n \| f(x_i)\| \;  \|f(c_i) - d_i\| \; \|
f(a_{n+1})\|\\
  & + & \sum_{i=1}^n \| f(x_i) - y_i \| (\| f(a_i) \| + \|
d_if(a_{n+1})\| )\\ 
& \leq & \sum_{i=1}^n (1 + \| y_i \|) \|f(c_i) - d_i \|  \; \|
f(a_{n+1})\|\\ 
 & + &  \sum_{i=1}^n \| f(x_i) - y_i \| (\| f(a_i) \| +
\delta(\varepsilon')/2)\\ 
& < & \varepsilon''n(1+N)\| f(a_{n+1}) \| + 
\varepsilon''n(M+\delta(\varepsilon')/2) \\
& \leq  & \delta(\varepsilon')\; .
\end{eqnarray*}   
Therefore  $\sum_{i=1}^n x_i(a_i + c_ia_{n+1})$ is invertible in $A$
(we use here  only the fact that $f$ is near-unit full) and thus
${\mathbf{a}} + {\mathbf{c}}a_{n+1}
\in Lg_n(A)$. We also have, for all $i = 1, \ldots, n$,  
\begin{eqnarray*} 
\| e - f(e + c_ia_{n+1}) \| & = & \| f(c_ia_{n+1}) \| \\
& \leq  & \|
f(c_ia_{n+1}) - d_if(a_{n+1}) \|  +  \| d_if(a_{n+1}) \| \\
 &\leq & \| f(c_i) - d_i \| \; \| f(a_{n+1}) \| +
\delta(\varepsilon')/2 \\
& \leq & \varepsilon'' \| f(a_{n+1}) \| + \delta(\varepsilon')/2 
\leq
\delta(\varepsilon')\; .  
\end{eqnarray*} 
Since $f$ is strong near-unit full, $e + c_ia_{n+1}$ is invertible 
in $A$ with an inverse $v_i \in A^{-1}$ and $\| v_i - e\| <
\varepsilon'$. Hence
$$\| c_ia_{n+1} \| = \| e - (e + c_ia_{n+1}) \| \leq \| v_i - e \| \|
e + c_ia_{n+1}
\| \leq \varepsilon'(1 + \| c_ia_{n+1} \|) .$$
This yields $\| c_ia_{n+1} \| \leq \varepsilon'/(1-\varepsilon') =
\varepsilon$. Therefore the conclusion  of Theorem 2.6 holds. The
same conclusion holds if $f(a_{n+1}) = 0$. Indeed, in this case 
$(f(a_1), \ldots, f(a_n)) \in
Lg_n(B)$ and, using the fact that $f$ is
near-unit full, $(a_1,
\ldots, a_n) \in Lg_n(A)$. This shows that
$$(a_1+z_1a_{n+1}, \ldots, 
a_n+z_na_{n+1}) \in Lg_n(A)  \mbox{ and }
\| z_i \| < \varepsilon \; ,$$  
for the choice $z_i = 0$, $i = 1, \ldots , n$. The proof is 
complete. 	\myqed 
\smallskip
 
The following result is a  
solution for Swan's problem for the absolute connected stable rank. 

\begin{Cor} 
Let $A$ and $B$ be two unital good normed algebras and 
$f : A \rightarrow B$ a strong near-unit full algebra morphism with
dense image. Then $acsr(A) = acsr(B)$. 
\end{Cor}

{\bf Proof :} It is known \cite[Appendix]{Ati67} that $f$ induces a
dense morphism $f : C(I,A) \rightarrow C(I,B)$. This new
morphism is also strong near-unit full. Indeed, let
$\delta(\varepsilon) > 0$ be such that
$\| f(a) - e \| < \delta(\varepsilon)$ implies $a$ invertible in $A$
and
$\| e - a^{-1} \| < \varepsilon$. Let now $g \in C(I,A)$ such that
$$\sup \{
\| f(g(t)) - e \| : t \in I \} < \delta (\varepsilon ) \; .$$ 
For every $t$
in $I$, the element $f(g(t))$ is $\delta(\varepsilon)$-close to the
unity of $B$ and thus $g(t)$ is invertible in $A$ with an inverse
$\varepsilon$-close to the unity of $A$. 
This shows that $g$ belongs to
$C(I,A^{-1}) = (C(I,A))^{-1}$ and its
inverse is $\varepsilon$-close to the unity of $C(I,A)$~; cf. the
proof of Lemma \ref{Lemma 3.2}. 

By Theorem \ref{snuf}, it follows that $tsr(C(I,A)) = tsr(C(I,B))$.
This, together with Theorem \ref{Nis}, completes the proof. \myqed 

\smallskip

We show now that Swan's problem for the connected stable rank has an affirmative 
answer at least for commutative Banach algebras.

\begin{Theo} 
Let $f : A \rightarrow B$ be an injective, dense and full
morphism  of commutative unital Banach algebras. Then $csr(A) =
csr(B)$.
\end{Theo}

{\bf Proof :} Let $X(A)$ and $X(B)$ be the maximal ideal spaces of
$A$ and $B$, respectively. We  have $csr(B) = \min \{ k : \mbox{
for all } n
\geq k, Lg_n(A) \mbox{ is connected }\}$. We have \cite{Tay76}
$\pi_0(Lg_n(B)) = [ X(B),{\bf C}^n_{\ast} ]$, the set of homotopy classes
of maps.  Here ${\bf C}^n_{\ast}$  denotes
${\bf C}^n \setminus \{ (0,\ldots,0)\}$. Since ${\mathbf{S}}^{2n-1}$
is a deformation  retract of ${\bf C}^n_{\ast}$, we have 
$
\pi_0(Lg_n(B)) = [ X(B) ,{\mathbf{S}}^{2n-1}]$.  
Since the monomorphism $f$ is dense and full, the adjoint operator
$f^{\ast} :  X(B) \rightarrow X(A)$ is bijective. Then
$$\pi_0(Lg_n(B)) = [ X(B) , {\mathbf{S}}^{2n-1}] = [ X(A) ,
{\mathbf{S}}^{2n-1}] =
\pi_0(Lg_n(A))\; .$$ This shows that $csr(A) = csr(B)$. 	\myqed


\section{The nonunital case} 
In the above sections, some stable ranks of unital
topological algebras were considered. We show briefly in this section how we can
extend these results to the non-unital case. Suppose that $A$ is non-unital and let
$A_1$ its unitization with unity $1$. Denote by 
$$Lg_n(A_1, A) = \{ {\mathbf{a}} \in Lg_n(A_1) : {\mathbf{a}} \equiv
{\mathbf{e}}_1
 \mbox{ mod } A^n \} \; ,$$  
where ${\mathbf{e}}_1 = (1, 0, \ldots, 0)$. Then (\cite{Rie83}, \cite{Pan90})
the Bass stable  rank $Bsr(A)$ of $A$ is defined as the least integer $n$ for which
the following condition holds: For every ${\mathbf{a}} = (a_1, \ldots, a_n)
\in Lg_{n+1}(A_1, A)$, there exists ${\mathbf{c}} = (c_1, \ldots, c_n) \in
A^n$  such that
$(a_1 + c_1a_{n+1}, \ldots, a_n + c_na_{n+1}) \in Lg_n(A_1, A)$. If no such $n$
exists, we set $Bsr(A) = \infty$ . The topological stable rank $tsr(A)$ is defined as
the lest integer $n$ such that $Lg_n(A_1, A)$ is dense in $\{ {\mathbf{a}}
\in A_1^n : {\mathbf{a}} \equiv {\mathbf{e}}_1 \mbox{ mod } A^n\}$. Then
\cite{Rie83}, \cite{Pan90}, $Bsr(A) = Bsr(A_1)$ and $tsr(A) = tsr(A_1)$. We
refer for instance to \cite{Pan90}, \cite{Vas69} for several considerations
of stable rank assertions for nonunital Banach algebras and rings without
identity element. 

Let $B$ be a Fr\'echet good algebra and $A$ a dense subalgebra. Let $\tilde{A}$ be
$A$ if $A$  has a unit and the unitization $A_1$ if not. If $A$ has a unit then $B$
has the same unit and we set $\tilde{B} = B$. For the nonunital case, we denote by 
$\tilde{B}$ the
algebra $B$ with the same unit adjoined (even if $B$ is unital already, adjoin a new
one). We say that a dense subalgebra $A$ of the algebra $B$ is
full if $a \in \tilde{A}$ is invertible in $\tilde{A}$ if it is invertible in
$\tilde{B}$. The following result is a consequence of the above results in the unital
case.

\begin{Theo}[Stability theorem] A dense and full 
Fr\'echet $\ast$-subalgebra of a $C^{\ast}$-algebra, closed under $C^{\infty}$
functional calculus of selfadjoint elements has the same Bass and topological stable
ranks. If the Fr\'echet subalgebra is a normed algebra, then they have also the same
absolute connected stable rank. For $tsr$ and $acsr$, A is considered with the
topology induced by that of $B$.
\end{Theo}


\end{document}